\documentclass{amsproc}

\usepackage{amsmath,amssymb}
\newtheorem{theorem}{Theorem}[section]
\newtheorem{lemm}[theorem]{Lemma}
\newtheorem{prop}[theorem]{Proposition}

\theoremstyle{definition}
\newtheorem{defi}[theorem]{Definition}

\newtheorem{coro}[theorem]{Corollary}
\theoremstyle{remark}
\newtheorem{remark}[theorem]{Remark}

\numberwithin{equation}{section}

\def\lg{\langle}
\def\rg{\rangle}

\def\al{\alpha}
\def\be{\beta}

\def\dim{\hbox{dim}}

\def\a{\alpha}

\newfont{\df}{eufm10}

\def\dim{\hbox{\rm dim}\,}

\begin{document}

\title[Quantizations of generalized-Witt and Jacobson-Witt algebras]
{Quantizations of generalized-Witt algebra and \\
of Jacobson-Witt algebra in the modular case}

\author[Hu]{Naihong Hu$^\star$}
\address{Department of Mathematics, East China Normal University,
Shanghai 200062, PR China} \email{nhhu@math.ecnu.edu.cn}

\author[Wang]{Xiuling Wang}
\address{Department of Mathematics,
East China Normal University,
Shanghai 200062, PR China\\
and Department of Mathematics, Harbin Normal University, Harbin
150080, PR China} \email{xlingwang@163.com}
\thanks{$^\star$N.H., Corresponding author,
supported by the NNSF (Grant 10431040), the PCSIRT, the TRAPOYT and
the FUDP from the MOE of China, and the SRSTP from the STCSM}

\date{Version on Feb. 12nd, 2006}

\dedicatory{To the memory of Professor Xihua Cao}
\begin{abstract}
We quantize the generalized-Witt algebra in characteristic $0$ with
its Lie bialgebra structures discovered by Song-Su (\cite{GY}). Via
a modulo $p$ reduction and a modulo ``$p$-restrictedness" reduction
process, we get $2^n{-}1$ families of truncated $p$-polynomial
noncocommutative deformations of the restricted universal enveloping
algebra of the Jacobson-Witt algebra $\mathbf{W}(n;\underline{1})$
(for the Cartan type simple modular restricted Lie algebra of $W$
type). They are new families of noncommutative and noncocommutative
Hopf algebras of dimension $p^{1+np^n}$ in characteristic $p$. Our
results generalize a work of Grunspan (J. Algebra 280 (2004),
145--161, \cite{CG}) in rank $n=1$ case in characteristic $0$. In
the modular case, the argument for a refined version follows from
the modular reduction approach (different from \cite{CG}) with some
techniques from the modular Lie algebra theory.
\end{abstract}

\maketitle
\section{Introduction and Definitions} 
\medskip

In a paper by Michaelis (\cite{M}) a class of infinite-dimensional
Lie bialgebras containing the Virasoro algebra was presented.
Afterwards, this type of Lie bialgebras was further classified by Ng
and Taft (\cite{NT}). Recently, Song and Su (\cite{GY}) classified
all Lie bialgebra structures on a given Lie algebra of generalized
Witt type (\cite{DZ}), which turned out to be coboundary triangular
(for the definition, see p. 28, \cite{ES}).

In Hopf algebra or quantum group theory, two standard methods to
yield new bialgebras from old ones are by twisting the product by a
$2$-cocycle but keeping the coproduct unchanged, and by twisting the
coproduct by a Drinfel'd twist but preserving the product.
Constructing quantizations of Lie bialgebras is an important
approach to producing new quantum groups (see \cite{D}, \cite{ES}
and references therein). Recently, Grunspan \cite{CG} obtained the
quantizations of the (infinite-dimensional) Witt algebra $W$ in
characteristic $0$ by using the twist discovered by Giaquinto and
Zhang \cite{AJ}, and of its simple modular Witt algebra
$\mathbf{W}(1;\underline{1})$ of dimension $p$ in characteristic $p$
by a reduction modulo $p$. In the modular case, however, his
treatment did not work for the restricted universal enveloping
algebra $\mathbf{u}(\mathbf{W}(1;\underline{1}))$. One reason is
that the periodic quotient with respect to the ideal $J$ he adopted
(see Definition 5, p. 158 \cite{CG}) resulted in some extra
relations like $[e_{p-2}, e_2]=4e_p\equiv 4e_o\ne 0$ incorrectly
imposed on it.
So in this case we need to consider the best way to approach it.
Another reason is that the assertion in Theorem 2 \cite{CG} (a
corrected version) is valid only for one specific case $i=1$ (except
for the other cases $i=2,\cdots,p{-}1$). This means that just one
family (rather than $p{-}1$ families asserted in \cite{CG}) of
polynomial noncocommutative deformations of
$\mathbf{u}(\mathbf{W}(1;\underline{1}))$ is able to be obtained via
a modular reduction process. The correct modular reduction approach
should be clarified as stated below.

In this paper, we first use the general quantization method by a
Drinfel'd twist (cf. \cite{CP}) to quantize explicitly the newly
defined triangular Lie bialgebra structure on the generalized-Witt
algebra in characteristic $0$ (\cite{GY}). Actually, this process
completely depends on the construction of Drinfel'd twists which, up
to integral scalars, are controlled by the classical Yang-Baxter
$r$-matrix. To study the modular case in characteristic $p$, we work
over the so-called {\it ``positive"} part subalgebra $\mathbf{W}^+$
of the generalized-Witt algebra $\mathbf{W}$. It is an
infinite-dimensional simple Lie algebra when defined over a field of
characteristic $0$, while, defined over a field of characteristic
$p$, it contains a maximal ideal $J_{\underline{1}}$ and the
corresponding quotient is exactly the Jacobson-Witt algebra
$\mathbf{W}(n;\underline{1})$, which is a Cartan type restricted
simple modular Lie algebra of $W$ type.  In order to yield the {\it
expected} finite-dimensional quantizations of the restricted
universal enveloping algebra for the Jacobson-Witt algebra
$\mathbf{W}(n;\underline{1})$, in concept at first, we need to give
a precise definition concerning what is the quantization of the
above object in the modular case (see Definition 3.3), and then go
on two steps of our reduction: {\it modulo $p$ reduction} and {\it
modulo ``$p$-restrictedness" reduction}. As a crucial immediate step
for the {\it modulo $p$ reduction}, we need first to develop the
quantization integral forms for the $\mathbb{Z}$-form
$\mathbf{W}_{\mathbb Z}^+$ in characteristic $0$. In this case, a
new phenomena appeared is that there exist $n$ the so-called
 {\it basic Drinfel'd twists} whose pairwise different products among them
afford the possible Drinfel'd twists (see Proposition 2.12).
Accordingly, we get $2^n{-}1$ new quantization integral forms for
$\mathbf{W}_{\mathbb Z}^+$ in characteristic $0$, which, via the
{\it modulo ``$p$-restrictedness" reduction}, eventually lead to
$2^n{-}1$ new examples of Hopf algebra of dimension $p^{1{+}np^n}$
with indeterminate $t$, or of dimension $p^{np^n}$ with specializing
$t$ into a scalar in $\mathcal{K}$ in characteristic $p$.

A remarkable point is that these Hopf algebras we obtained contain
the well-known Radford algebra as a Hopf subalgebra. Our work
extends the class of examples of noncommutative and noncocommutative
finite-dimensional Hopf algebras in finite characteristic (see
\cite{T}).

\subsection{Generalized-Witt algebra and its Lie bialgebra structure}
Let $\mathbb{F}$ be a field with $\text{char}(\mathbb{F})=0$ and let
$n>0$. The notations used here are the same as in \cite{DZ}. Let
$\mathbb{F}[x^{\pm1}_1,\cdots,x^{\pm1}_n]$ be a Laurent polynomial
algebra, and $\partial_i$ coincides with the degree operator
$x_i\frac{\partial}{\partial x_i}$ .
 Set
$T=\bigoplus_{i=1}^n{\mathbb{Z}}\partial_i$. For
$\alpha=(\alpha_1,\cdots,\alpha_n)
 \in \mathbb{Z}^{n}$, write $x^{\alpha}=x_1^{\alpha_1} \cdots x_n^{\alpha_n}$.

Denote
$\mathbf{W}=\mathbb{F}[x^{\pm1}_1,\cdots,x^{\pm1}_n]\otimes_{\mathbb{Z}}
T=\text{Span}_{\mathbb{F}}\{\,x^{\alpha}\partial \mid \alpha \in
\mathbb{Z}^{n},\partial\in T\}$, where we set $x^{\alpha}\partial\
=x^{\alpha}\otimes
\partial$ for short. Then $\mathbf{W}=\text{Der}_{\mathbb{F}}(\mathbb{F}[x^{\pm1}_1,\cdots,x^{\pm1}_n])$ is a Lie algebra of
generalized-Witt type (for definition, see \cite{DZ}) under the
following bracket
$$[\,x^{\alpha}\partial,\,x^{\beta}\partial'\,]
 =x^{\alpha+\beta}\bigl(\,\partial(\beta)\partial'-\partial'(\alpha)\partial\,\bigr),
 \qquad\forall\
 \alpha,\, \beta\in \mathbb{Z}^{n};\ \partial,\, \partial' \in T,$$
where
$\partial(\beta)=\lg\partial,\beta\rg=\lg\beta,\partial\rg=\sum\limits_{i=1}^{n}a_i\beta_i\in\mathbb{Z}$
for $\partial=\sum\limits_{i=1}^{n}a_i\partial_i \in T$ and
$\beta=(\beta_1,\cdots,\beta_n) \in \mathbb{Z}^{n}$. The bilinear
map $\lg\cdot,\cdot\,\rg: \,T\times {\mathbb{Z}}^n\longrightarrow
{\mathbb{Z}}$ is non-degenerate in the sense:
$\partial(\al)=\lg\partial,\al\rg=0$ $(\forall\;\partial \in T),\
\Longrightarrow \al=0,$
and $\partial(\al)=\lg\partial,\al\rg=0$  $(\forall\;\al \in
\mathbb{Z}^n),\ \Longrightarrow
\partial=0.$
$\mathbf{W}$ is an infinite dimensional simple Lie algebra over
$\mathbb{F}$.

\smallskip
The following result is due to \cite{GY}.
\begin{prop}
There is a triangular Lie bialgebra structure on $\mathbf{W}$
given by the classical Yang-Baxter $r$-matrix
$r:=\partial_0\otimes x^{\gamma}\partial'_0-x^{\gamma}\partial'_0
\otimes \partial_0$, for $\partial_0, \;\partial'_0\in T, \;\gamma
\in \mathbb{Z}^{n}$, where $[\,\partial_0,
x^\gamma\partial'_0\,]=\partial_0(\gamma)x^\gamma\partial_0'$.
\end{prop}

From Proposition 1.1, we notice that the classical Yang-Baxter
$r$-matrix $r$ is uniquely expressed as the antisymmetric tensor of
{\it two distinguished elements} $\partial_0,\, x^\gamma\partial'_0$
up to scalars satisfying $[\,\partial_0,
x^\gamma\partial'_0\,]=\partial_0(\gamma)x^\gamma\partial'_0$. In
fact, for a given $r$-matrix, we may take {\it two distinguished
elements} of the form $h:=\partial_0(\gamma)^{-1}\partial_0$ and
$e:=\partial_0(\gamma)x^{\gamma}\partial'_0$ such that $[h, e]=e$,
where $\partial_0(\gamma)\neq 0\;(\in\mathbb{Z})$.

\subsection{Generalized-Witt subalgebra $\mathbf{W}^+$}
Denote $D_i=\frac{\partial}{\partial x_i}$. Set
$\mathbf{W}^+:=\text{Span}_{\mathcal{K}}\{x^{\alpha} D_i\mid
\alpha\in\mathbb{Z}_+^n, 1\le i\le n\}$, where $\mathbb{Z}_+$ is
the set of non-negative integers. Then
$\mathbf{W}^+=\text{Der}_{\mathcal{K}}(\mathcal{K}[x_1,\cdots,x_n])$
is the derivation Lie algebra of polynomial ring
$\mathcal{K}[x_1,\cdots,x_n]$, which, via the identification
$x^\al D_i$ with $x^{\alpha-\epsilon_i}
\partial_i$, can be considered as a Lie
subalgebra (the ``{\it positive}" part) of the generalized-Witt
algebra $\mathbf{W}$ over a field $\mathcal{K}$. Evidently, we have
the following result (see Ex. 8, p. 153 in \cite{HR})
\begin{lemm} $(1)$ If $\mathcal{K}=\mathbb{F}$, i.e., $\text{char}(\mathcal{K})=0$, then
$\mathbf{W}^+$ is a simple Lie algebra of infinite dimension,
which is a Lie subalgebra of $\mathbf{W}$.

$(2)$ If $\text{char}(\mathcal{K})=p$, then there exists a maximal
ideal $J_{\underline{1}}:=\langle \{x^\alpha D_i\mid \exists j:
\alpha_j\ge p, 1\le i\le n\}\rangle$ in $\mathbf{W}^+$ such that
$\mathbf{W}^+/J_{\underline{1}}\cong \mathbf{W}(n; \underline{1})$
under the identification $\frac{1}{\alpha!}x^\al D_i$ $(0\le
\al\le\tau)$ with $x^{(\al)}D_i$ and the others with $0$, where
$\mathbf{W}(n; \underline{1})$ is the Jacobson-Witt algebra defined
in the subsection below.
\end{lemm}

\subsection{Jacobson-Witt algebra $\mathbf{W}(n;\underline{1})$}

Assume now that $\text{char}(\mathcal{K})=p$, then by definition
(cf. \cite{H}), the Jacobson-Witt algebra
$\mathbf{W}(n;\underline{1})$ is a restricted simple Lie algebra
over a field $\mathcal{K}$. Its structure of $p$-Lie algebra is
given by $D^{[p]}=D^p,\; \forall\, D \in
\mathbf{W}(n;\underline{1})$ with a basis $\{\,x^{(\alpha)}D_j \mid
1\leq j\leq n, \ 0 \leq \alpha \leq \tau \}$, where
$\underline{1}=(1,\cdots,1),\, \tau=(p{-}1,\cdots,p{-}1) \in
\mathbb{N}^n$; $\epsilon_i=(\delta_{1i},\cdots,\delta_{ni})$ such
that $x^{(\epsilon_i)}=x_i$ and $D_j(x_i)=\delta_{ij}$; and
$\mathcal{O}(n;\underline{1}):=\{\,x^{(\al)}\mid 0 \leq \alpha \leq
\tau \}$ is the restricted divided power algebra with
$x^{(\al)}x^{(\be)}=\binom{\al{+}\be}{\al}\,x^{(\al{+}\be)}$ and a
convention: $x^{(\alpha)}=0$ if $\alpha$ has a component
$\alpha_j<0$ or $>p$, where
$\binom{\al{+}\be}{\al}:=\prod_{i=1}^n\binom{\al_i{+}\be_i}{\al_i}$.
Moreover, $\mathbf{W}(n;\underline{1})\cong
\text{Der}_{\mathcal{K}}(\mathcal{O}(n;\underline{1}))$.

By definition, the restricted universal enveloping algebra
$\mathbf{u}(\mathbf{W}(n;\underline{1}))$ is isomorphic to
$U(\mathbf{W}(n;\underline{1}))/I$ where $I$ is the Hopf ideal of
$U(\mathbf{W}(n;\underline{1}))$ generated by $H_i^p-H_i,\, D^p$
with $D\neq H_i=x^{(\epsilon_i)}D_i$. Since
$dim_{\mathcal{K}}\mathbf{W}(n;\underline{1})=np^n$, we have
$dim_{\mathcal{K}}\mathbf{u}(\mathbf{W}(n;\underline{1}))$
$=p^{np^n}$.

\subsection{Quantization by Drinfel'd twists} The following
result is well-known (see \cite{CP}, \cite{D}, \cite{ES}, etc.).
\begin{lemm}
Let $(A,m,\iota,\Delta_0,\varepsilon_0,S_0)$ be a Hopf algebra over
a commutative ring. A Drinfel'd twist $\mathcal{F}$ on $A$ is an
invertible element of $A\otimes A$ such that
\begin{gather*}(\mathcal{F}\otimes
1)(\Delta_0\otimes \text{\rm Id})(\mathcal{F})=(1\otimes
\mathcal{F})(\text{\rm Id}\otimes\Delta_0)(\mathcal{F}), \\
(\varepsilon_0\otimes \text{\rm Id})(\mathcal{F})=1=(\text{\rm
Id}\otimes \varepsilon_0)(\mathcal{F}).
\end{gather*} Then,
$w=m(\text{\rm Id}\otimes S_0)(\mathcal{F})$ is invertible in $A$
with $w^{-1}=m(S_0\otimes \text{\rm Id})(\mathcal{F}^{-1})$.

Moreover, if we define $\Delta: \,A\longrightarrow A\otimes A$ and
$S: \,A\longrightarrow A$ by
$$\Delta(a)=\mathcal{F}\Delta_0(a)\mathcal{F}^{-1},
\qquad S=w\,S_0(a)\,w^{-1},$$ then $(A, m, \iota,
\Delta,\varepsilon,S)$ is a new Hopf algebra, called the twisting
of $A$ by the Drinfel'd twist $\mathcal{F}$.
\end{lemm}

Let $\mathbb{F}[[t]]$ be a ring of formal power series over a field
$\mathbb{F}$ with $\text{char}(\mathbb{F})=0$. Assume that $L$ is a
triangular Lie bialgebra over $\mathbb{F}$ with a classical
Yang-Baxter $r$-matrix $r$ (see \cite{D}, \cite{ES}). Let $U(L)$ be
the universal enveloping algebra of $L$, with the standard Hopf
algebra structure $(U(L),m,\iota,\Delta_0,\varepsilon_0,S_0)$.

\smallskip
Now let us consider {\it the topologically free
$\mathbb{F}[[t]]$-algebra} $U(L)[[t]]$ (for definition, see p. 4,
\cite{ES}), which can be viewed as an associative
$\mathbb{F}$-algebra of formal power series with coefficients in
$U(L)$. Naturally, $U(L)[[t]]$ equips with an induced Hopf algebra
structure arising from that on $U(L)$. By abuse of notation, we
denote it by $(U(L)[[t]],m,\iota,\Delta_0,\varepsilon_0,S_0)$.

\begin{defi}
For a triangular Lie bialgebra $L$ over $\mathbb{F}$ with
$\text{char}(\mathbb{F})=0$, $U(L)[[t]]$ is called {\it a
quantization of $U(L)$ by a Drinfel'd twist} $\mathcal{F}$ over
$U(L)[[t]]$ if $U(L)[[t]]/tU(L)[[t]]\cong U(L)$, and $\mathcal{F}$
is determined by its $r$-matrix $r$ (namely, its Lie bialgebra
structure).
\end{defi}

\subsection{A crucial Lemma}
For any element $x$ of a unital $R$-algebra ($R$ a ring) and
     $a \in R$, we set (see \cite{AJ})
$$x_a^{\lg n\rg}:=(x+a)(x+a+1)\cdots(x+a+n-1) \leqno(1)$$
and $x^{\lg n\rg}:=x_0^{\lg n\rg}$.

We also set
$$x_a^{[n]}:=(x+a)(x+a-1)\cdots(x+a-n+1) \leqno(2)$$
and $x^{[n]}:=x_0^{[n]}$.

\begin{lemm} $($\cite{AJ}, \cite{CG}$)$
For any element $x$ of a unital $\mathbb{F}$-algebra with
$\text{char}(\mathbb{F})=0$, $a, \,b \in \mathbb{F}$ and $r,\,
s,\, t \in \mathbb{Z}$, one has
\begin{gather*}
x_a^{\lg s+t\rg}=x_a^{\lg s\rg}\,x_{a+s}^{\lg t\rg},\tag{3} \\
x_a^{[s+t]}=x_a^{[s]}\,x_{a-s}^{[t]},\tag{4} \\
x_a^{[s]}=x_{a-s+1}^{\lg s\rg},\tag{5} \\
\sum\limits_{s+t=r}\frac{(-1)^t}{s!\,t!}\,x_a^{[s]}\,x_b^{\lg
t\rg}=\dbinom{a{-}b} {r}=\frac{(a{-}b)\cdots(a{-}b{-}r{+}1)}{r!}, \tag{6} \\
\sum\limits_{s+t=r}\frac{(-1)^t}{s!\,t!}\,x_a^{[s]}\,x_{b-s}^{[t]}=\dbinom{a{-}b{+}r{-}1}
{r}=\frac{(a{-}b)\cdots(a{-}b{+}r{-}1)}{r!}.\tag{7}
\end{gather*}
\end{lemm}

\section{Quantization of Lie bialgebra of generalized-Witt type}
\subsection{Some commutative relations in $U(\mathbf{W})$}
For the universal enveloping algebra $U(\mathbf{W})$ of the
generalized-Witt algebra $\mathbf{W}$ over $\mathbb{F}$, we need to
do some necessary calculations, which are important for the
construction of the Drinfel'd twists in the sequel.

\begin{lemm} Fix the two distinguished
elements $h:=\partial_0(\gamma)^{-1}\partial_0$ and
$e:=\partial_0(\gamma)x^{\gamma}\partial'_0$ with
$\partial_0(\gamma)\in\mathbb{Z}^*$ for $\mathbf{W}$. For $a \in
\mathbb{F},\,\alpha \in \mathbb{Z}^{n}$, and $m,\,k$ non-negative
integers, the following equalities hold in $U(\mathbf{W}):$
\begin{gather*}
x^{\alpha}\partial\cdot
h_a^{[m]}=h_{a-\frac{\partial_0(\alpha)}{\partial_0(\gamma)}}^{[m]}\cdot
x^{\alpha}\partial,
\tag{8}\\
x^{\alpha}\partial\cdot h_a^{\lg
m\rg}=h_{a-\frac{\partial_0(\alpha)}{\partial_0(\gamma)}}^{\lg
m\rg}\cdot x^{\alpha}\partial, \tag{9}\\
e^k\cdot h_a^{[m]}=h_{a-k}^{[m]}\cdot e^k, \tag{10}\\
e^k\cdot h_a^{\lg m\rg}=h_{a-k}^{\lg m\rg}\cdot e^k, \tag{11}\\
x^{\alpha}\partial\cdot(x^{\beta}\partial')^m=\sum\limits_{\ell=0}^m({-}1)^\ell\dbinom{m}
{\ell} (x^{\beta}\partial')^{m{-}\ell}\cdot
x^{\alpha{+}\ell\beta}\Bigl(a_\ell
\partial-b_\ell\partial'\Bigr), \tag{12}
\end{gather*}
where
$a_\ell=\prod\limits_{j=0}^{\ell-1}\partial'(\alpha{+}j\beta)$,
$b_\ell=\ell\,\partial(\beta)a_{\ell{-}1}$, and set $a_0=1$,
$b_0=0$.
\end{lemm}
\begin{proof}
One has $x^{\alpha}\partial\cdot\partial_0=\partial_0\cdot
x^{\alpha}\partial -\partial_0(\alpha)\,x^{\alpha}\partial$. So it
is easy to see that (8) is true for $m=1$. Suppose that (8) is
true for $m$, then
\begin{equation*}
\begin{split}
x^{\alpha}\partial\cdot h_a^{[m+1]} &= x^{\alpha}\partial \cdot
h^{[m]}_a\cdot h_{a-m}\\ & =
           h^{[m]}_{a-\frac{\partial_0(\alpha)}{\partial_0(\gamma)}}\cdot x^{\alpha}\partial \cdot
           (h+a-m)\\
&=
h^{[m]}_{a-\frac{\partial_0(\alpha)}{\partial_0(\gamma)}}\cdot\Bigl(\partial_0
(\gamma)^{-1}(\partial_0\cdot x^{\alpha}\partial
-\partial_0(\alpha)x^{\alpha}\partial)+(a-m)x^{\alpha}\partial\Bigr)
\\
 &=
 h^{[m+1]}_{a-\frac{\partial_0(\alpha)}{\partial_0(\gamma)}}\cdot x^{\alpha}\partial.
\end{split}
\end{equation*}
Hence (8) holds for all $m$. Similarly, we can get (9), (10) and
(11) by induction.

Formula (12) is a consequence of the fact (see Proposition 1.3
(4), \cite{HR}) that for any elements $a,\,c$ in an associative
algebra, one has
$$c\,a^m=\sum_{\ell=0}^m(-1)^{\ell}\dbinom{m}{\ell}a^{m{-}\ell}(\text{ad}\,a)^\ell(c),$$
together with the formula
$$
(\text{ad}\,x^{\beta}\partial')^\ell(x^\alpha\partial)
=x^{\alpha{+}\ell\beta}(a_\ell\partial-b_\ell\partial'),\leqno(13)$$
obtained by induction on $\ell$ when taking $a=x^{\beta}\partial',
\,c=x^\alpha\partial$.
\end{proof}

Let us denote by $(U(\mathbf{W}),\, m,\, \iota,\, \Delta_0,\,
S_0,\, \varepsilon_0)$ the standard Hopf algebra structure on
$U(\mathbf{W})$, $i.e.$, we have the definitions of the coproduct,
the antipode and the counit as follows
\begin{gather*}
\Delta_0(x^{\alpha}\partial)=x^{\alpha}\partial\otimes 1+1\otimes
 x^{\alpha}\partial,\\
 S_0(x^{\alpha}\partial)=-x^{\alpha}\partial,\\
 \varepsilon_0(x^{\alpha}\partial)=0.
\end{gather*}

\subsection{Quantization of $U(\mathbf{W})$ in characteristic 0}
To describe a quantization of $U(\mathbf{W})$ by a Drinfel'd twist
$\mathcal{F}$ over $U(\mathbf{W})[[t]]$, we need to construct
explicitly such a Drinfel'd twist. In what follows, we shall see
that such a twist in our case depends heavily upon the choice of
{\it two distinguished elements} $h,\,e$ arising from its
$r$-matrix $r$ (see subsection 1.1).

Now we proceed with the construction. For $a \in \mathbb{F}$, we set
\begin{gather*}
\mathcal{F}_a=\sum\limits_{r=0}^{\infty}\frac{(-1)^r}{r!}h_a^{[r]}\otimes
e^rt^r,\qquad F_a=\sum\limits_{r=0}^{\infty}\frac{1}{r!}h_a^{\lg
r\rg}\otimes
e^rt^r,\\
u_a=m\cdot(S_0\otimes \text{\rm Id})(F_a),\qquad\qquad
v_a=m\cdot(\text{\rm Id}\otimes S_0)(\mathcal{F}_a).
\end{gather*}

Write $\mathcal{F}=\mathcal{F}_0,\, F=F_0,\,u=u_0,\,v=v_0$.

Since $S_0(h_a^{\lg r\rg})=(-1)^rh_{-a}^{[r]}$ and
$S_0(e^r)=(-1)^re^r$, we obtain
$$v_a=\sum\limits_{r=0}^{\infty}\frac{1}{r!}h_a^{[r]}
e^rt^r, \quad
u_b=\sum\limits_{r=0}^{\infty}\frac{(-1)^r}{r!}h_{-b}^{[r]}
e^rt^r.$$
\begin{lemm}
For $a,b \in \mathbb{F},$ one has $$ \mathcal{F}_a
F_b=1\otimes(1-et)^{a-b}, \quad\text{and }\quad v_a
u_b=(1-et)^{-(a+b)}.$$
\end{lemm}
\begin{proof}
Using Lemma 1.5, we obtain
\begin{equation*}
\begin{split}
\mathcal{F}_a F_b &=
\sum\limits_{r,s=0}^{\infty}\frac{(-1)^r}{r!s!}h_a^{[r]}h_b^{\lg
s\rg}\otimes e^re^st^rt^s
\\ &=\sum\limits_{m=0}^{\infty}(-1)^m\left(\sum\limits_{r+s=m}\frac{(-1)^s}{r!s!}h_a^{[r]}h_b^{\lg s\rg}\right)\otimes
e^mt^m \\
&=\sum\limits_{m=0}^{\infty}(-1)^m\dbinom{a-b}{m}\otimes e^mt^m
=1\otimes (1-et)^{a-b}.
\\
v_a u_b &=
\sum\limits_{m,n=0}^{\infty}\frac{(-1)^m}{m!n!}h_a^{[n]}e^nh_{-b}^{[m]}e^mt^{m+n}
\\ &=\sum\limits_{r=0}^{\infty}\sum\limits_{m+n=r}\frac{(-1)^m}{m!n!}h_a^{[n]}h_{-b-n}^{[m]}
e^rt^r \\
&=\sum\limits_{r=0}^{\infty}\dbinom{a+b+r-1}{r} e^rt^r
=(1-et)^{-(a+b)}.
\end{split}
\end{equation*}
This completes the proof.
\end{proof}
\begin{coro}
For $a \in \mathbb{F}$, $\mathcal{F}_a$ and $u_a$ are invertible
with $\mathcal{F}_a^{-1}=F_a$ and $u_a^{-1}=v_{-a}$.  In particular,
$\mathcal{F}^{-1}=F$ and $u^{-1}=v$.
\end{coro}
\begin{lemm}\ For any positive integers $r$, we have
$$\Delta_0(h^{[r]})=\sum\limits_{i=0}^r \dbinom{r}{i}h^{[i]}\otimes
h^{[r-i]}.$$  Furthermore, $\Delta_0(h^{[r]})=\sum\limits_{i=0}^r
\dbinom{r}{i}h^{[i]}_{-s}\otimes h^{[r-i]}_s$ for any $s \in
\mathbb{F}$.
\end{lemm}
\begin{proof}
By induction on $r$, it is easy to see that it is true for $r=1$.
If it is true for $r$, then
\begin{equation*}
\begin{split}
\Delta_0(h^{[r+1]}) &=\left(\sum\limits_{i=0}^r
\dbinom{r}{i}h^{[i]}\otimes h^{[r-i]}\right)\Bigl((h-r)\otimes
1+1\otimes (h-r)+r(1\otimes 1)\Bigr)
\\
&=\left(\sum\limits_{i=1}^{r-1} \dbinom{r}{i}h^{[i]}\otimes
h^{[r-i]}\right)((h-r)\otimes 1+1\otimes (h-r)) \\
& \quad+r\left(
\sum\limits_{i=0}^r \dbinom{r}{i}h^{[i]}\otimes
h^{[r-i]}\right)+1\otimes h^{[r+1]}
+(h-r)\otimes h^{[r]} \\
& \quad +h^{[r+1]}\otimes 1+h^{[r]}\otimes(h-r) \\
&=1\otimes h^{[r+1]}+h^{[r+1]}\otimes 1+ r\left(
\sum\limits_{i=1}^{r-1} \dbinom{r}{i}h^{[i]}\otimes
h^{[r-i]}\right) \\
& \quad +h\otimes h^{[r]}+h^{[r]}\otimes h+\sum\limits_{i=1}^{r-1}
\dbinom{r}{i}h^{[i+1]}\otimes h^{[r-i]}\\
& \quad + \sum\limits_{i=1}^{r-1}(i-r) \dbinom{r}{i}h^{[i]}\otimes
h^{[r-i]}+\sum\limits_{i=1}^{r-1} \dbinom{r}{i}h^{[i]}\otimes
h^{[r-i+1]}\\
&\quad+\sum\limits_{i=1}^{r-1}(-i) \dbinom{r}{i}h^{[i]}\otimes
h^{[r-i]}
\\ &=1\otimes h^{[r+1]}+h^{[r+1]}\otimes 1+\sum\limits_{i=1}^{r}
\left[\dbinom{r}{i{-}1}+\dbinom{r}{i}\right]h^{[i]}\otimes h^{[r-i+1]} \\
&= \sum\limits_{i=0}^{r+1}
\dbinom{r+1}{i}h^{[i]}\otimes h^{[r+1-i]}.
\end{split}
\end{equation*}
Therefore, the formula holds by induction.
\end{proof}
\begin{prop}
$\mathcal{F}=\sum\limits_{r=0}^{\infty}\frac{(-1)^r}{r!}h^{[r]}\otimes
e^rt^r$ is a Drinfel'd twist on $U(\mathbf{W})[[t]]$, that is to
say, the equalities hold
\begin{gather*}
(\mathcal{F}\otimes 1)(\Delta_0\otimes \text{\rm
Id})(\mathcal{F})=(1\otimes
\mathcal{F})(\text{\rm Id}\otimes\Delta_0)(\mathcal{F}),\\
(\varepsilon_0\otimes \text{\rm Id})(\mathcal{F})=1=(\text{\rm
Id}\otimes \varepsilon_0)(\mathcal{F}).
\end{gather*}
\end{prop}
\begin{proof} The second equality holds evidently. For the first one, by Lemmas 2.4 \& 1.5, it is easy to get
\begin{equation*}
\begin{split}
\text{LHS} &=(\mathcal{F}\otimes 1)(\Delta_0\otimes \text{\rm
Id})(\mathcal{F})
 \\ &=
\left(\sum\limits_{r=0}^{\infty}\frac{(-1)^rt^r}{r!}h^{[r]}\otimes
e^r\otimes 1\right)
\left(\sum\limits_{s=0}^{\infty}\frac{(-1)^st^s}{s!}\sum\limits_{i=0}^s
\dbinom{s}{i}h^{[i]}_{-r}\otimes h_r^{[s-i]}\otimes e^s \right)
\\ &=\sum\limits_{r,s=0}^{\infty}\frac{(-1)^{r+s}t^{r+s}}{r!s!}\sum\limits_{i=0}^s
\dbinom{s}{i}h^{[r]}h^{[i]}_{-r}\otimes e^rh_r^{[s-i]}\otimes e^s
\\ &=\sum\limits_{r,s=0}^{\infty}\frac{(-1)^{r+s}t^{r+s}}{r!s!}\sum\limits_{i=0}^s
\dbinom{s}{i}h^{[r+i]}\otimes h^{[s-i]}e^r\otimes e^s.
\\
\end{split}
\end{equation*}
\begin{equation*}
\begin{split}
\text{RHS} &=(1\otimes \mathcal{F})(\text{\rm
Id}\otimes\Delta_0)(\mathcal{F})
\\ &=
\left(\sum\limits_{r=0}^{\infty}\frac{(-1)^rt^r}{r!} 1\otimes
h^{[r]}\otimes e^r\right)
\left(\sum\limits_{s=0}^{\infty}\frac{(-1)^st^s}{s!}\sum\limits_{i=0}^s
\dbinom{s}{i}h^{[s]}\otimes e^i\otimes e^{s-i} \right) \\
&=\sum\limits_{r,s=0}^{\infty}\frac{(-1)^{r+s}t^{r+s}}{r!s!}\sum\limits_{i=0}^s
\dbinom{s}{i}h^{[s]}\otimes h^{[r]}e^i\otimes e^{r+s-i}.
\end{split}
\end{equation*}
It suffices to show that
$$\sum\limits_{p+q=m}\frac{1}{p!q!}\sum\limits_{j=0}^p
\dbinom{p}{j}h^{[q+j]}\otimes h^{[p-j]}e^q \otimes e^p
=\sum\limits_{r+s=m}\frac{1}{r!s!}\sum\limits_{i=0}^s
\dbinom{s}{i}h^{[s]}\otimes h^{[r]}e^i\otimes e^{m-i}.
$$
Fix $r, \,s,\, i$ such that $r+s=m,\, 0\leq i \leq s$. Set
$q=i,\,q+j=s$, then $p=m-i,\,p-j=r$. We can see that the
coefficients of $h^{[s]}\otimes h^{[r]}e^i\otimes e^{m-i}$ in both
sides are equal.
\end{proof}
Having the above Proposition in hand, by Lemma 1.3, now we can
perform the process of twisting the standard Hopf structure
$(U(\mathbf{W})[[t]],\, m,\, \iota,\, \Delta_0,\, S_0,\,
\varepsilon_0)$ by the Drinfel'd twist $\mathcal{F}$ constructed
above.

\smallskip
The following Lemma is very useful to our main result in this
section.
\begin{lemm}
For $a \in \mathbb{F},\,\alpha \in \mathbb{Z}^n$, one has
\begin{gather*}
((x^{\alpha}\partial)^s\otimes 1)\cdot
F_a=F_{a-s\frac{\partial_0(\alpha)}{\partial_0(\gamma)}}\cdot
((x^{\alpha}\partial)^s\otimes 1), \tag{14}\\
(1\otimes x^{\alpha}\partial)\cdot
F_a=\sum\limits_{\ell=0}^{\infty}(-1)^\ell F_{a+\ell}\cdot
\Bigl(h_a^{\lg \ell\rg}\otimes
x^{\alpha+\ell\gamma}(A_\ell\partial-B_\ell\partial'_0)t^\ell\Bigr), \tag{15}\\
(x^{\alpha}\partial)\cdot
u_a=u_{a+\frac{\partial_0(\alpha)}{\partial_0(\gamma)}}\cdot
\Bigl(\sum\limits_{\ell=0}^{\infty}x^{\alpha+\ell\gamma}(A_\ell\partial
-B_\ell\partial'_0)\cdot h_{1-a}^{\lg \ell\rg}t^\ell\Bigr), \tag{16}\\
(x^{\alpha}\partial)^s\cdot
u_a=u_{a+s\frac{\partial_0(\alpha)}{\partial_0(\gamma)}}\cdot
\Bigl(\sum\limits_{\ell=0}^{\infty}d^{(\ell)}\bigl((x^\alpha\partial)^s\bigr)\cdot
h_{1-a}^{\lg \ell\rg}t^\ell\Bigr), \tag{17} \\
(1\otimes (x^{\alpha}\partial)^s)\cdot
F_a=\sum\limits_{\ell=0}^{\infty}(-1)^\ell F_{a+\ell}\cdot
\Bigl(h_a^{\lg \ell\rg}\otimes
d^{(\ell)}((x^{\alpha}\partial)^s)t^\ell\Bigr), \tag{18}
\end{gather*}
where $d^{(\ell)}:=\frac{1}{\ell!}(\text{\rm ad}\,e)^{\ell}$,
$A_\ell=\frac{\partial_0(\gamma)^\ell}{\ell!}
\prod\limits_{j=0}^{\ell{-}1}\partial'_0(\alpha{+}j\gamma),
\,B_\ell=\partial_0(\gamma)\partial(\gamma)A_{\ell{-}1}$, and set
$A_0=1$, $B_0=0$.
\end{lemm}
\begin{proof}
For (14): By (9), one has
\begin{equation*}
\begin{split}
(x^\alpha\partial\otimes 1)\cdot F_a
&=\sum\limits_{m=0}^{\infty}\frac{1}{m!}x^\alpha\partial \cdot
h_a^{\lg m\rg}\otimes e^mt^m
\\
&=\sum\limits_{m=0}^{\infty}\frac{1}{m!}h_{a-\frac{\partial_0(\alpha)}{\partial_0{(\gamma)}}}^{\lg
m\rg}\cdot x^\alpha\partial \otimes e^mt^m
\\
&=F_{a-\frac{\partial_0(\alpha)}{\partial_0{(\gamma)}}}\cdot(x^\alpha\partial
\otimes 1).
\end{split}
\end{equation*}
By induction on $s$, we obtain the result.

\noindent For (15): Let
$a_\ell=\prod\limits_{j=0}^{\ell-1}\partial'_0(\alpha{+}j\gamma),\,
b_\ell=\ell\,\partial(\gamma)a_{\ell{-}1}$. Then using (12) we get
\begin{equation*}
\begin{split}
(1 \otimes x^\alpha\partial )\cdot F_a
&=\sum\limits_{m=0}^{\infty}\frac{1}{m!} h_a^{\lg m\rg}\otimes
x^\alpha\partial \cdot e^mt^m
\\
&=\sum\limits_{m=0}^{\infty}\frac{1}{m!}h_{a}^{\lg m\rg}\otimes
\left(
\sum\limits_{\ell=0}^{m}(-1)^\ell\dbinom{m}{\ell}\partial_0(\gamma)^\ell
e^{m{-}\ell}\cdot
x^{\alpha{+}\ell\gamma}(a_\ell\partial{-}b_\ell\partial'_0)t^m\right)
\\
&=\sum\limits_{m=0}^{\infty}\sum\limits_{\ell=0}^{\infty}(-1)^\ell\frac{1}{m!\ell!}h_{a}^{\lg
m+\ell\rg}\otimes
\partial_0(\gamma)^\ell e^{m}\cdot x^{\alpha+\ell\gamma}(a_\ell\partial-b_\ell\partial'_0)t^{m+\ell}
\\
&=\sum\limits_{\ell=0}^{\infty}(-1)^\ell\left(\sum\limits_{m=0}^{\infty}\frac{1}{m!}h_{a+\ell}^{\lg
m\rg}\otimes e^mt^m\right)\Bigl(h_a^{\lg \ell\rg}\otimes
x^{\alpha+\ell\gamma}(A_\ell\partial-B_\ell\partial'_0)t^{\ell}\Bigr)
\\ &=\sum\limits_{\ell=0}^{\infty}(-1)^\ell F_{a+\ell}\cdot \Bigl(h_a^{\lg \ell\rg}\otimes
x^{\alpha+\ell\gamma}(A_\ell\partial-B_\ell\partial'_0)t^{\ell}\Bigr).
\end{split}
\end{equation*}
For (16): Using (8), (10) and (12), we get
\begin{equation*}
\begin{split}
x^{\alpha}\partial\cdot u_a &=x^{\alpha}\partial\cdot
\left(\sum\limits_{r=0}^{\infty}\frac{(-1)^r}{r!} h_{-a}^{[r]}\cdot
e^rt^r \right) \\ &=
\sum\limits_{r=0}^{\infty}\frac{(-1)^r}{r!}x^{\alpha}\partial\cdot
h_{-a}^{[r]}\cdot e^rt^r  \\
&=\sum\limits_{r=0}^{\infty}\frac{(-1)^r}{r!}h_{-a-\frac{\partial_0(\alpha)}
{\partial_0{(\gamma)}}}^{[r]}\cdot x^{\alpha}\partial\cdot
 e^rt^r \\ &=\sum\limits_{r=0}^{\infty}\frac{(-1)^r}{r!}h_{-a-\frac{\partial_0(\alpha)}{\partial_0{(\gamma)}}}^{[r]}
 \left(\sum\limits_{\ell=0}^{r}(-1)^\ell\dbinom{r}{\ell}\partial_0(\gamma)^\ell e^{r-\ell}\cdot
 x^{\alpha+\ell\gamma}(a_\ell\partial-b_\ell\partial'_0)t^r
 \right)
 \\ &=\sum\limits_{r,\ell=0}^{\infty}\frac{({-}1)^{r{+}\ell}}{(r{+}\ell)!}
 h_{{-}a{-}\frac{\partial_0(\alpha)}{\partial_0{(\gamma)}}}^{[r{+}\ell]}
 \left(({-}1)^\ell\dbinom{r{+}\ell}{\ell}\partial_0(\gamma)^\ell e^{r}\cdot
 x^{\alpha{+}\ell\gamma}(a_\ell\partial{-}b_\ell\partial'_0)t^{r{+}\ell}
  \right)\\
&=\sum\limits_{r,\ell=0}^{\infty}\frac{(-1)^{r}}{r!\ell!}\partial_0(\gamma)^\ell
  h_{-a-\frac{\partial_0(\alpha)}
  {\partial_0{(\gamma)}}}^{[r]}\cdot h_{-a-\frac{\partial_0(\alpha)}
  {\partial_0{(\gamma)}}-r}^{[\ell]}\cdot
 e^{r}\cdot x^{\alpha+\ell\gamma}(a_\ell\partial-b_\ell\partial'_0)t^{r+\ell}
  \\
&=\sum\limits_{\ell=0}^{\infty}\left(\sum\limits_{r=0}^{\infty}\frac{(-1)^{r}}{r!}
h_{-a-\frac{\partial_0(\alpha)}
  {\partial_0{(\gamma)}}}^{[r]}\cdot e^{r}t^r\right)\cdot h_{-a-\frac{\partial_0(\alpha)}
  {\partial_0{(\gamma)}}}^{[\ell]}\cdot x^{\alpha+\ell\gamma}(A_\ell\partial-B_\ell\partial'_0)t^{\ell}
  \\
&=u_{a+\frac{\partial_0(\alpha)}
  {\partial_0{(\gamma)}}}\cdot\sum\limits_{\ell=0}^{\infty}h_{-a-\frac{\partial_0(\alpha)}
  {\partial_0{(\gamma)}}}^{[\ell]}\cdot x^{\alpha+\ell\gamma}(A_\ell\partial-B_\ell\partial'_0)t^{\ell}
  \\
&=u_{a+\frac{\partial_0(\alpha)}
  {\partial_0{(\gamma)}}}\cdot\sum\limits_{\ell=0}^{\infty}x^{\alpha+\ell\gamma}
  (A_\ell\partial-B_\ell\partial'_0)\cdot h_{-a+\ell}^{[\ell]}t^{\ell}
 \\
 &=u_{a+\frac{\partial_0(\alpha)}
  {\partial_0{(\gamma)}}}\cdot\sum\limits_{\ell=0}^{\infty}x^{\alpha+\ell\gamma}
  (A_\ell\partial-B_\ell\partial'_0)\cdot h_{-a+1}^{\lg \ell\rg}t^{\ell}.
\end{split}
\end{equation*}
As for (17): By (13), we get
$$d^{(\ell)}(x^\alpha\partial)=x^{\alpha+\ell\gamma}(A_\ell\partial-B_\ell\partial_0').\leqno(19)$$
Combining with (16), we get that (17) is true for $s=1$. Using the
derivation property of $d^{(\ell)}$, we easily obtain that
$$d^{(\ell)}(a_1\cdots
a_s)=\sum_{\ell_1{+}\cdots{+}\ell_s=\ell}d^{(\ell_1)}(a_1)\cdots
d^{(\ell_s)}(a_s).\leqno(20)$$
By induction on $s$, we
have
\begin{equation*}
\begin{split}
(x^{\alpha}\partial)^{s{+}1}\cdot u_a &= x^{\alpha}\partial\cdot
u_{a+ s\frac{\partial_0(\alpha)}{\partial_0 (\gamma)}} \cdot
\sum\limits_{n=0}^{\infty}d^{(n)}((x^{\alpha}\partial)^s)\cdot
h_{1-a}^{\langle
n \rangle} t^{n}\\
&= u_{a{+}(s{+}1)\frac{\partial_0(\alpha)}{\partial_0
(\gamma)}}{\cdot}\Bigl(\sum\limits_{m=0}^{\infty}d^{(m)}(x^{\alpha}\partial)\cdot
h_{1{-}a{-}s\frac{\partial_0(\alpha)}{\partial_0
(\gamma)}}^{\langle m \rangle}
t^{m}\Bigr)\\
&\quad\qquad\qquad\qquad\cdot\Bigl(\sum\limits_{n=0}^{\infty}d^{(n)}((x^{\alpha}\partial)^s)\cdot
h_{1{-}a}^{\langle
n \rangle} t^{n}\Bigr)\\
&=u_{a{+}(s{+}1)\frac{\partial_0(\alpha)}{\partial_0
(\gamma)}}\cdot\Bigl(\sum\limits_{m,n=0}^{\infty}d^{(m)}(x^{\alpha}\partial)
d^{(n)}((x^{\alpha}\partial)^s)h_{1-a+n}^{\langle m
\rangle} h_{1-a}^{\langle n \rangle} t^{n+m}\Bigr)\\
&=u_{a{+}(s{+}1)\frac{\partial_0(\alpha)}{\partial_0
(\gamma)}}\cdot\Bigl(\sum\limits_{\ell=0}^{\infty}\sum\limits_{m+n=\ell}d^{(m)}(x^{\alpha}\partial)
d^{(n)}((x^{\alpha}\partial)^s) h_{1-a}^{\langle \ell \rangle} t^{\ell}\Bigr)\\
&=u_{a{+}(s{+}1)\frac{\partial_0(\alpha)}{\partial_0
(\gamma)}}\cdot\Bigl(\sum\limits_{\ell=0}^{\infty}d^{(\ell)}((x^{\alpha}\partial)^{s{+}1})h_{1-a}^{\langle
\ell \rangle} t^{\ell}\Bigr),
\end{split}
\end{equation*}
where we get the first and second  ``=" by using the inductive
hypothesis and (16), the third by using (13) \& (20) and the
fourth by using (4) \& (20).

For (18): this follows from (15) \& (20).

Thus,  the proof is complete.
\end{proof}

The following theorem gives the quantization of $U(\mathbf{W})$ by
the Drinfel'd twist $\mathcal{F}$, which is essentially determined
by the Lie bialgebra triangular structure on $\mathbf{W}$.
\begin{theorem}
With the choice of two distinguished elements
$h:=\partial_0(\gamma)^{-1}\partial_0$,
$e:=\partial_0(\gamma)x^{\gamma}\partial'_0$
$(\partial_0(\gamma)\in \mathbb{Z}^*)$ with $[h,\,e]=e$
 in the generalized-Witt
algebra $\mathbf{W}$ over $\mathbb{F}$, there exists a structure
of noncommutative and noncocommutative Hopf algebra
$(U(\mathbf{W})[[t]],m,\iota,\Delta,S,\varepsilon)$ on
$U(\mathbf{W})[[t]]$ over $\mathbb{F}[[t]]$
 with
$U(\mathbf{W})[[t]]/tU(\mathbf{W})[[t]]$ $\cong U(\mathbf{W})$,
which leaves the product of $U(\mathbf{W})[[t]]$ undeformed but
with the deformed coproduct, antipode and counit defined by
\begin{gather*}
\Delta(x^{\alpha}\partial)=x^{\alpha}\partial\otimes
(1-et)^{\frac{\partial_0(\alpha)}{\partial_0(\gamma)}}+\sum\limits_{\ell=0}^{\infty}(-1)^\ell
h^{\lg \ell\rg}\otimes (1-et)^{-\ell}\cdot
x^{\alpha+\ell\gamma}(A_\ell\partial-B_\ell\partial'_0)t^\ell,
\tag{21}\\
S(x^{\alpha}\partial)=-(1-et)^{-\frac{\partial_0(\alpha)}{\partial_0(\gamma)}}\cdot\Bigl(\sum\limits_{\ell=0}^{\infty}
x^{\alpha+\ell\gamma}(A_\ell\partial-B_\ell\partial'_0)\cdot
h_1^{\lg \ell\rg}t^\ell\Bigr),
\tag{22}\\
\varepsilon(x^{\alpha}\partial)=0,  \tag{23}
\end{gather*}
where $\alpha \in \mathbb{Z}^n,\,
A_\ell=\frac{\partial_0(\gamma)^\ell}{\ell!}\prod\limits_{j=0}^{\ell-1}\partial'_0(\alpha{+}j\gamma),\,
B_\ell=\partial_0(\gamma)\partial(\gamma)A_{\ell{-}1}$, and set
$A_0=1$, $B_0=0$.
\end{theorem}
\begin{proof}
By Lemmas 1.3 and 2.2, it follows from (14) and (15) that
\begin{equation*}
\begin{split}
\Delta(x^{\alpha}\partial)
&=\mathcal{F}\cdot\Delta_0(x^{\alpha}\partial)\cdot\mathcal{F}^{-1}\\
&=\mathcal{F}\cdot(x^{\alpha}\partial\otimes 1)\cdot F+
\mathcal{F}\cdot(1\otimes x^{\alpha}\partial)\cdot F
\\
&=\Bigl(\mathcal{F}
F_{-\frac{\partial_0(\alpha)}{\partial_0(\gamma)}}\Bigr)\cdot(x^{\alpha}\partial{\otimes}
1)+ \sum\limits_{\ell=0}^{\infty}(-1)^\ell
\Bigl(\mathcal{F}F_{\ell}\Bigr)\cdot\Bigl(h^{\lg \ell\rg}{\otimes}
x^{\alpha{+}\ell\gamma}(A_\ell\partial{-}B_\ell\partial'_0)t^\ell\Bigr)
\\ &=\Bigl(1\otimes
(1{-}et)^{\frac{\partial_0(\alpha)}{\partial_0(\gamma)}}\Bigr)\cdot
(x^{\alpha}\partial\otimes
1)\\
&\quad+\sum\limits_{\ell=0}^{\infty}(-1)^\ell \Bigl(1\otimes
(1{-}et)^{-\ell}\Bigr)\cdot\Bigl( h^{\lg \ell\rg}\otimes
x^{\alpha{+}\ell\gamma}(A_\ell\partial{-}B_\ell\partial'_0)t^\ell\Bigr)
\\ &=x^{\alpha}\partial\otimes
(1{-}et)^{\frac{\partial_0(\alpha)}{\partial_0(\gamma)}}+\sum\limits_{\ell=0}^{\infty}(-1)^\ell
h^{\lg \ell\rg}\otimes (1{-}et)^{-\ell}\cdot
x^{\alpha+\ell\gamma}(A_\ell\partial{-}B_\ell\partial'_0)t^\ell.
\\
\end{split}
\end{equation*}
By (16) and Lemma 2.2, we obtain
\begin{equation*}
\begin{split}
S(x^{\alpha}\partial)&=u^{-1}S_0(x^{\alpha}\partial)\,u=-v\cdot
x^{\alpha}\partial\cdot u\\ &=-v \cdot
u_{\frac{\partial_0(\alpha)}{\partial_0(\gamma)}}\cdot
\Bigl(\sum\limits_{\ell=0}^{\infty}x^{\alpha+\ell\gamma}(A_\ell\partial
-B_\ell\partial'_0)\cdot h_{1}^{\lg \ell\rg}t^\ell\Bigr)
\\ &=-(1-et)^{-\frac{\partial_0(\alpha)}{\partial_0(\gamma)}}\cdot\Bigl(\sum\limits_{\ell=0}^{\infty}
x^{\alpha+\ell\gamma}(A_\ell\partial-B_\ell\partial'_0)\cdot
h_1^{\lg \ell\rg}t^\ell\Bigr).
\end{split}
\end{equation*}
Hence, we get the result.
\end{proof}

For later use, we need to make the following
\begin{lemm} For $s\ge 1$, one has
\begin{gather*}
\Delta((x^\alpha\partial)^s)=\sum_{0\le j\le s\atop
\ell\ge0}\dbinom{s}{j}({-}1)^\ell(x^{\alpha}\partial)^jh^{\langle
\ell\rangle}\otimes(1{-}et)^{j\frac{\partial_0(\al)}{\partial_0(\gamma)}{-}\ell}
d^{(\ell)}((x^\al\partial)^{s{-}j})t^\ell.\tag{\text{\rm i}}\\
S((x^{\alpha}\partial)^s)=
(-1)^s(1-et)^{-s\frac{\partial_0(\alpha)}{\partial_0(\gamma)}}\cdot\Bigl(\sum\limits_{\ell=0}^{\infty}
d^{(\ell)}((x^{\alpha}\partial)^s)\cdot h_1^{\lg
\ell\rg}t^\ell\Bigr).\tag{\text{\rm ii}}
\end{gather*}
\end{lemm}
\begin{proof} By (18) and Lemma 2.2, we obtain
\begin{equation*}
\begin{split}
\Delta((x^{\alpha}\partial)^s)&=\mathcal{F}\Bigl(x^{\al}\partial\otimes
1+1\otimes x^{\alpha}\partial\Bigr)^s\mathcal{F}^{-1}\\
&=\sum_{j=0}^s\binom{s}{j}\mathcal{F}F_{-j\frac{\partial_0(\alpha)}{\partial_0(\gamma)}}
(x^\alpha\partial{\otimes}
1)^j\Bigl(\sum_{\ell\ge0}({-}1)^\ell\mathcal{F}F_{\ell}\Bigl(h^{\langle\ell\rangle}{\otimes}
d^{(\ell)}((x^\alpha\partial)^{s{-}j})t^\ell\Bigr)\Bigr)\\
&=\sum_{j=0}^s\sum_{\ell\ge0}\binom{s}{j}({-}1)^\ell\Bigl((x^\alpha\partial)^j{\otimes}
(1{-}et)^{j\frac{\partial_0(\al)}{\partial_0(\gamma)}{-}\ell}\Bigr)\Bigl(h^{\langle\ell\rangle}{\otimes}
d^{(\ell)}((x^\alpha\partial)^{s{-}j})t^\ell\Bigr)\\
&=\sum_{0\le j\le s\atop
\ell\ge0}\dbinom{s}{j}({-}1)^\ell(x^{\alpha}\partial)^jh^{\langle
\ell\rangle}\otimes(1{-}et)^{j\frac{\partial_0(\al)}{\partial_0(\gamma)}{-}\ell}
d^{(\ell)}((x^\al\partial)^{s{-}j})t^\ell.
\end{split}
\end{equation*}
Again by (17) and Lemma 2.2, we get
\begin{equation*}
\begin{split}
S((x^{\alpha}\partial)^s)&=u^{-1}S_0((x^{\alpha}\partial)^s)\,u=(-1)^s
v\cdot (x^{\alpha}\partial)^s\cdot u\\ &=(-1)^s v \cdot
u_{s\frac{\partial_0(\alpha)}{\partial_0(\gamma)}}\cdot
\Bigl(\sum\limits_{\ell=0}^{\infty}d^{(\ell)}((x^{\alpha}\partial)^s)\cdot
h_{1}^{\lg \ell\rg}t^\ell\Bigr)
\\ &=(-1)^s(1-et)^{-s\frac{\partial_0(\alpha)}{\partial_0(\gamma)}}\cdot\Bigl(\sum\limits_{\ell=0}^{\infty}
d^{(\ell)}((x^{\alpha}\partial)^s)\cdot h_1^{\lg
\ell\rg}t^\ell\Bigr).
\end{split}
\end{equation*}
This  completes the proof.
\end{proof}

\subsection{Quantization integral forms of $\mathbb{Z}$-form $\mathbf{W}_{\mathbb{Z}}^+$ in
characteristic $0$} Note that $\{x^\alpha D_i\mid
\alpha\in\mathbb{Z}_+^n\}$ is a $\mathbb{Z}$-basis of
$\mathbf{W}_{\mathbb{Z}}^+$. In order to get the quantization
integral forms of $\mathbb{Z}$-form $\mathbf{W}_{\mathbb{Z}}^+$, it
suffices to consider what conditions are for those coefficients
occurred in the formulae (21) \& (22) to be integral for the
indicated basis elements.

Write $r=\partial_0(\gamma)$. Notice that suitable powers of
factors $(1-et)^{\pm\frac{1}{r}}$ and coefficients $A_\ell$ and
$B_\ell$ occur in  (21) \& (22). Grunspan (\cite{CG}) proved
\begin{lemm}
For any $a,\,k,\,\ell\in\mathbb{Z}$,
$a^\ell\prod\limits_{j=0}^{\ell-1}(k{+}ja)/\ell!$ is an integer.
\end{lemm}
According to this Lemma, we see that, if we take
$\partial'_0=\partial_0$, then $A_\ell$ and $B_\ell$ are integers.
However, those coefficients occurred in the expansions of
$(1-et)^{\pm\frac{1}{r}}$ are not integral unless $r=1$.
Consequently, the case we are interested in is only when
$h=\partial_k$, $e=x^{\epsilon_k}\partial_k$ ($1\le k\le n$),
namely, $r=1$. Denote by $\mathcal{F}(k)$ the corresponding
Drinfel'd twist. As a result of Theorem 2.7, we have
\begin{coro}
With the specific choice of the two distinguished elements
$h:=x^{\epsilon_k}D_k$, $e:=x^{2\epsilon_k}D_k$ $(1\le k\le n)$ with
$[h,\,e]=e$, the corresponding quantization integral form of
$U(\mathbf{W}^+_{\mathbb{Z}})$ over
$U(\mathbf{W}^+_{\mathbb{Z}})[[t]]$ by the Drinfel'd twist
$\mathcal{F}(k)$ with the algebra structure undeformed is given by
\begin{gather*}
\Delta(x^{\alpha}D_i)=x^{\alpha}D_i{\otimes}
(1{-}et)^{\alpha_k{-}\delta_{ik}}{+}\sum\limits_{\ell=0}^{\infty}{({-}1)}^\ell
C_\ell\,h^{\lg \ell\rg}{\otimes} (1{-}et)^{{-}\ell}\cdot
x^{\alpha{+}\ell\epsilon_k}D_it^\ell,
\tag{24}\\
S(x^{\alpha}D_i)={-}(1{-}et)^{-\alpha_k{+}\delta_{ik}}\cdot\Bigl(\sum\limits_{\ell=0}^{\infty}
C_\ell\,x^{\alpha{+}\ell\epsilon_k}D_i\cdot h_1^{\lg
\ell\rg}t^\ell\Bigr),
\tag{25}\\
\varepsilon(x^{\alpha}D_i)=0,  \tag{26}
\end{gather*}
where $\alpha \in \mathbb{Z}^n_+,\,C_\ell=A_\ell-B_\ell,\,
A_\ell=\frac{1}{\ell!}\prod\limits_{j=0}^{\ell-1}(\alpha_k{-}\delta_{ik}{+}j),\,
B_\ell=\delta_{ik}A_{\ell{-}1}$, and set $A_0=1$, $B_0=0$.
\end{coro}
\begin{remark}
It is interesting to consider those products
$\mathcal{F}(j_1)\cdots\mathcal{F}(j_s)$ of some pairwise different
Drinfel'd twists $\mathcal{F}(j_1),\cdots,\mathcal{F}(j_s)$ with
$j_1<\cdots<j_s$ in the system of the $n$ so-called {\it basic
Drinfel'd twists} $\{\,\mathcal{F}(1),\cdots,\mathcal{F}(n)\,\}$.
Using the same argument as in the proof of Theorem 2.7, one can get
many more new Drinfel'd twists (which depend on a bit more
calculations to be carried out), which will not only lead to many
more new quantization integral forms over the
$U(\mathbf{W}_{\mathbb{Z}}^+)[[t]]$, but the possible quantizations
over the $\mathbf{u}_{t,q}(\mathbf{W}(n;\underline{1}))$ as well,
via our modular reduction approach developed in the next section.
\end{remark}
More precisely, we note that $[\mathcal{F}(i), \mathcal{F}(j)]=0$
for any $1\le i, j\le n$. This fact, according to the definition of
$\mathcal{F}(i)$, implies the commutative relations in the case
$i\ne j$:
\begin{equation*}
\begin{split}
(\mathcal{F}(j)\otimes 1)(\Delta_0\otimes\text{\rm Id})
(\mathcal{F}(i))&=(\Delta_0\otimes\text{\rm Id})
(\mathcal{F}(i))(\mathcal{F}(j)\otimes 1),\\
(1\otimes \mathcal{F}(j))(\text{\rm Id}\otimes\Delta_0)
(\mathcal{F}(i))&=(\text{\rm Id}\otimes\Delta_0)
(\mathcal{F}(i))(1\otimes\mathcal{F}(j)),
\end{split}\tag{\text{$*$}}
\end{equation*}
which give rise to the following property.
\begin{prop}
$\mathcal{F}(i)\mathcal{F}(j)(i\neq j)$ is still a Drinfel'd twist
on $U(\mathbf{W}_{\mathbb{Z}}^+)[[t]]$. In general,
$\mathcal{F}^\eta:=\mathcal{F}(1)^{\eta_1}\cdots\mathcal{F}(n)^{\eta_n}(\eta_i=0$
or $1)$ is a Drinfel'd twist on $U(\mathbf{W}_{\mathbb{Z}}^+)[[t]]$.
\end{prop}
\begin{proof}
Note that $\Delta_0\otimes\text{\rm id}$, $\text{\rm
id}\otimes\Delta_0$, $\varepsilon_0\otimes\text{\rm id}$ and
$\text{\rm id}\otimes\varepsilon_0$ are algebraic homomorphisms.
According to Proposition 2.5, it suffices to show that
$$
(\mathcal{F}(i)\mathcal{F}(j)\otimes 1)(\Delta_0\otimes \text{\rm
Id})(\mathcal{F}(i)\mathcal{F}(j)) =(1\otimes
\mathcal{F}(i)\mathcal{F}(j))(\text{\rm Id}\otimes\Delta_0)
(\mathcal{F}(i)\mathcal{F}(j)).
$$
Using $(*)$, we have
\begin{equation*}
\begin{split}
\text{\rm LHS}&=(\mathcal{F}(i)\otimes 1)(\mathcal{F}(j)\otimes
1)(\Delta_0\otimes \text{\rm Id})(\mathcal{F}(i))(\Delta_0\otimes
\text{\rm
Id})(\mathcal{F}(j))\\
&=(\mathcal{F}(i)\otimes 1)(\Delta_0\otimes \text{\rm
Id})(\mathcal{F}(i))(\mathcal{F}(j)\otimes 1)(\Delta_0\otimes
\text{\rm
Id})(\mathcal{F}(j))\\
&=(1\otimes\mathcal{F}(i) )(\text{\rm Id}\otimes\Delta_0
)(\mathcal{F}(i))(1\otimes \mathcal{F}(j))(\text{\rm
Id}\otimes\Delta_0)(\mathcal{F}(j))\\
&=(1\otimes\mathcal{F}(i) )(1\otimes \mathcal{F}(j))(\text{\rm
Id}\otimes\Delta_0 )(\mathcal{F}(i))(\text{\rm
Id}\otimes\Delta_0)(\mathcal{F}(j))=\text{\rm RHS}.
\end{split}
\end{equation*}
As a result, $\mathcal{F}^\eta$ is also a Drinfel'd twist, which
gives a quantization of $U(\mathbf{W}_{\mathbb{Z}}^+)$.
\end{proof}
\begin{lemm}
With the specific choice of the distinguished elements
$h(k):=x^{\epsilon_k}D_k$, $e(k):=x^{2\epsilon_k}D_k$ with
$[h(k),\,e(k)]=e(k)$ and $h(m):=x^{\epsilon_m}D_m$,
$e(m):=x^{2\epsilon_m}D_m$ $(1\le m\neq k\le n)$ with
$[h(m),\,e(m)]=e(m)$, the corresponding quantization integral form
of $U(\mathbf{W}^+_{\mathbb{Z}})$ over
$U(\mathbf{W}^+_{\mathbb{Z}})[[t]]$ by the Drinfel'd twist
$\mathcal{F}=\mathcal{F}(m)\mathcal{F}(k)$ with the algebra
structure undeformed is given by
\begin{gather*}
\Delta(x^{\alpha}D_i)=x^{\alpha}D_i{\otimes}
(1{-}e(k)t)^{\alpha_k{-}\delta_{ik}}(1{-}e(m)t)^{\alpha_m{-}\delta_{im}}
{+}\sum\limits_{\ell,n=0}^{\infty}{({-}1)}^{\ell+n} \tag{27} \\
\cdot C(k)_\ell C(m)_n h(k)^{\lg \ell\rg} h(m)^{\lg n \rg} {\otimes}
(1{-}e(k)t)^{{-}\ell} (1{-}e(m)t)^{{-}n}\cdot
x^{\alpha{+}\ell\epsilon_k+n\epsilon_m}D_it^{\ell+n},\\
S(x^{\alpha}D_i)={-}(1{-}e(k)t)^{-\alpha_k{+}\delta_{ik}}(1{-}e(m)t)^{-\alpha_m{+}\delta_{im}}
 \Bigl(\sum\limits_{\ell,n=0}^{\infty} C(k)_\ell C(m)_n \tag{28} \\ \cdot
x^{\alpha{+}\ell\epsilon_k+n\epsilon_m}D_i\cdot h(m)_1^{\lg
n\rg}h(k)_1^{\lg \ell\rg}t^{\ell+n}\Bigr),\\
\varepsilon(x^{\alpha}D_i)=0,  \tag{29}
\end{gather*}
where $\alpha \in \mathbb{Z}^n_+,\,C(k)_\ell=A(k)_\ell
{-}B(k)_\ell,\, A(k)_\ell,\, B(k)_\ell$ as in Corollary 2.10.
\end{lemm}
\begin{proof} Using Corollary 2.10, we get
\begin{equation*}
\begin{split}
\Delta(x^{\alpha}D_i)&=\mathcal{F}(m)\mathcal{F}(k)
\Delta_0(x^{\alpha}D_i)\mathcal{F}(k)^{-1}\mathcal{F}(m)^{-1}\\
&=\mathcal{F}(m)\Bigl( x^{\alpha}D_i{\otimes}
(1{-}e(k)t)^{\alpha_k{-}\delta_{ik}}\\
&\qquad+\sum\limits_{\ell=0}^{\infty}{({-}1)}^\ell C(k)_\ell
h(k)^{\lg \ell\rg}{\otimes} (1{-}e(k)t)^{{-}\ell}\cdot
x^{\alpha{+}\ell\epsilon_k}D_it^\ell\Bigr)  \mathcal{F}(m)^{-1}
\end{split}
\end{equation*}
 Using (14) and Lemma 2.2, we get
\begin{equation*}
\begin{split}
\mathcal{F}(m)&\Bigl( x^{\alpha}D_i{\otimes}
(1{-}e(k)t)^{\alpha_k{-}\delta_{ik}}\Bigr) \mathcal{F}(m)^{-1}\\
&=\mathcal{F}(m)\Bigl( x^{\alpha}D_i{\otimes}1 \Bigr)
\mathcal{F}(m)^{-1}
\Bigl(1\otimes(1{-}e(k)t)^{\alpha_k{-}\delta_{ik}}\Bigr)
 \\&=\mathcal{F}(m)\mathcal{F}(m)^{-1}_{\delta_{im}{-}\alpha_m}
\Bigl( x^{\alpha}D_i{\otimes}1 \Bigr)\Bigl(1\otimes(1{-}e(k)t)^{\alpha_k{-}\delta_{ik}}\Bigr)\\
  &=\Bigl(1\otimes(1{-}e(m)t)^{\alpha_m{-}\delta_{im}}\Bigr)
 \Bigl( x^{\alpha}D_i{\otimes}1 \Bigr)
 \Bigl(1\otimes(1{-}e(k)t)^{\alpha_k{-}\delta_{ik}}\Bigr)\\
&=x^{\alpha}D_i\otimes
(1{-}e(k)t)^{\alpha_k{-}\delta_{ik}}(1{-}e(m)t)^{\alpha_m{-}\delta_{im}}.
\end{split}
\end{equation*}
Using (15), we get
\begin{equation*}
\begin{split}
\mathcal{F}(m)&\Bigl(\sum\limits_{\ell=0}^{\infty}{({-}1)}^\ell
C(k)_\ell h(k)^{\lg \ell\rg}{\otimes}(1{-}e(k)t)^{{-}\ell}\cdot
x^{\alpha{+}\ell\epsilon_k}D_it^\ell\Bigr) \mathcal{F}(m)^{-1}\\
&=\sum\limits_{\ell=0}^{\infty}{({-}1)}^\ell C(k)_\ell(h(k)^{\lg
\ell\rg}{\otimes}(1{-}e(k)t)^{{-}\ell})\cdot
\mathcal{F}(m)\Bigl(1\otimes
x^{\alpha{+}\ell\epsilon_k}D_it^\ell\Bigr) \mathcal{F}(m)^{-1}\\
&=\sum\limits_{\ell,n=0}^{\infty}{({-}1)}^{\ell+n}
C(k)_\ell(h(k)^{\lg \ell\rg}{\otimes} (1{-}e(k)t)^{{-}\ell})
\\& \quad  \cdot \mathcal{F}(m)\mathcal{F}(m)_{n}^{-1}\cdot
\Bigl(C(m)_nh(m)^{\lg n \rg}\otimes
x^{\alpha+\ell\epsilon_k+n\epsilon_m} D_it^{n+\ell}\Bigr)
\\&=\sum\limits_{\ell,n=0}^{\infty}{({-}1)}^{\ell+n}C(k)_\ell C(m)_n h(k)^{\lg \ell\rg}
h(m)^{\lg n \rg}{\otimes}  \\
&\qquad (1{-}e(k)t)^{{-}\ell}(1{-}e(m)t)^{{-}n}\cdot
x^{\alpha{+}\ell\epsilon_k+n\epsilon_m}D_it^{\ell+n}.
\end{split}
\end{equation*}
For $k\neq m$, using the definitions of $v$ and  $u$, we get
$v=v(k)v(m)=v(m)v(k)$ and $u=u(m)u(k)=u(k)u(m)$.
By
Corollary 2.10 and using (16), we have
\begin{equation*}
\begin{split}
S(x^{\alpha}D_i)&=-v\cdot x^{\alpha}D_i\cdot u =-v(m)v(k)\cdot x^{\alpha}D_i\cdot  u(k)u(m)\\
&=v(m)\cdot\Bigl({-}(1{-}e(k)t)^{-\alpha_k{+}\delta_{ik}}\cdot\Bigl(\sum\limits_{\ell=0}^{\infty}
C(k)_\ell x^{\alpha{+}\ell\epsilon_k}D_i\cdot h(k)_1^{\lg
\ell\rg}t^\ell\Bigr)\Bigr)\cdot u(m)\\
&={-}(1{-}e(k)t)^{-\alpha_k{+}\delta_{ik}}\cdot
v(m)u(m)_{\alpha_m-\delta_{im}}
\\
&\quad \cdot \Bigl(\sum\limits_{\ell, n=0}^{\infty}C(k)_{\ell}C(m)_n
x^{\alpha+\ell\epsilon_k+n\epsilon_m}D_i\cdot h(m)_{1}^{\lg
n\rg}h(k)_{1}^{\lg \ell \rg}t^{n+\ell}\Bigr)\\
&={-}(1{-}e(k)t)^{-\alpha_k{+}\delta_{ik}}(1{-}e(m)t)^{-\alpha_m{+}\delta_{im}}\\
&\quad \cdot \Bigl(\sum\limits_{\ell, n=0}^{\infty}C(k)_{\ell}C(m)_n
x^{\alpha+\ell\epsilon_k+n\epsilon_m}D_i\cdot h(m)_{1}^{\lg
n\rg}h(k)_{1}^{\lg \ell \rg}t^{n+\ell}\Bigr).
\end{split}
\end{equation*}
Therefore, the proof is complete.
\end{proof}

Let $I=\{1,\cdots,n\}$, $\eta=(\eta_1,\cdots,\eta_n)$ ($\eta_i
\in\{0, 1\}$), $\ell=(\ell_1,\cdots,\ell_n) \in \mathbb{N}^n$. Then
$|\eta\ell|=\sum_i\eta_i\ell_i$. Set
$(1{-}e(I)t)^{\eta(\alpha{-}\delta_{iI})}=\prod\limits_{k \in
I}(1{-}e(k)t)^{\eta_k(\alpha_k{-}\delta_{ik})}$ and $h(I)^{\lg
\eta\ell\rg}=\prod\limits_{k \in I}h(k)^{\lg \eta_k\ell_k\rg}$.

More generally, we can get the following result.
\begin{theorem} The
corresponding quantization integral form of
$U(\mathbf{W}^+_{\mathbb{Z}})$ over
$U(\mathbf{W}^+_{\mathbb{Z}})[[t]]$ by the Drinfel'd twist
$\mathcal{F}^{\eta}=\mathcal{F}(1)^{\eta_1}\cdots
\mathcal{F}(n)^{\eta_n}$ with the algebra structure undeformed is
given by
\begin{gather*}
\Delta(x^{\alpha}D_i)=x^{\alpha}D_i{\otimes}
(1{-}e(I)t)^{\eta(\alpha{-}\delta_{iI})}
{+}\sum\limits_{\eta\ell\ge\underline{0}}{({-}1)}^{|\eta\ell|}
h(I)^{\lg \eta\ell\rg}{\otimes}\tag{30} \\
\quad C(I)^{\eta}_\ell (1{-}e(I)t)^{{-}\eta\ell}\cdot
x^{\alpha{+}\eta\ell}D_it^{|\eta\ell|},
\\
S(x^{\alpha}D_i)={-}(1{-}e(I)t)^{-\eta(\alpha{-}\delta_{iI})}
\Bigl(\sum\limits_{\eta\ell\ge\underline{0}} C(I)^{\eta }_\ell
x^{\alpha{+}\eta\ell}D_i\cdot h(I)_1^{\lg
\eta\ell\rg}t^{|\eta\ell|}\Bigr), \tag{31} \\
\varepsilon(x^{\alpha}D_i)=0,  \tag{32}
\end{gather*}
where $\alpha \in \mathbb{Z}^n_+,$\, $C(I)_\ell=\prod\limits_{k \in
I}(A(k)_{\ell_k} {-}B(k)_{\ell_k})^{\eta_k}$,\, $A(k)_{\ell_k},\,
B(k)_{\ell_k}$ as in Corollary 2.10.
\end{theorem}

\medskip
\section{Quantization of the Jacobson-Witt algebra in the modular case}
\medskip

In this section, our main purposes are twofold: Firstly, in view of
Lemma 1.2, we make {\it the modulo $p$ reduction} for the
quantization integral form of $U(\mathbf{W}^+_{\mathbb{Z}})$ in
characteristic $0$ obtained in Corollary 2.10 to yield the
quantization of $U(\mathbf{W}(n;\underline{1}))$, for the Cartan
type restricted simple modular Lie algebra
$\mathbf{W}(n;\underline{1})$ of $W$ type in characteristic $p$.
Secondly, we shall further make {\it the ``$p$-restrictedness"
reduction} for the quantization of $U(\mathbf{W}(n;\underline{1}))$,
which will lead to the required quantization of
$\mathbf{u}(\mathbf{W}(n;\underline{1}))$ (here
$\mathbf{u}(\mathbf{W}(n;\underline{1}))$ is the restricted
universal enveloping algebra of $\mathbf{W}(n;\underline{1})$).

\subsection{Modulo $p\,$ reduction}
Let $\mathbb{Z}_p$ be the prime subfield of $\mathcal{K}$ with
$\text{char}(\mathcal{K})=p$. When considering
$\mathbf{W}_{\mathbb{Z}}^+$ as a $\mathbb{Z}_p$-Lie algebra, namely,
making a modulo $p$ reduction for the defining relations of
$\mathbf{W}_{\mathbb{Z}}^+$, we denote it by
$\mathbf{W}_{\mathbb{Z}_p}^+$. By Lemma 1.2 (2), we see that
$(J_{\underline{1}})_{\mathbb{Z}_p}=\text{Span}_{\mathbb{Z}_p}\{x^\alpha
D_i\mid \exists j: \alpha_j\ge p\,\}$ is a maximal ideal of
$\mathbf{W}^+_{\mathbb{Z}_p}$, and
$\mathbf{W}^+_{\mathbb{Z}_p}/(J_{\underline{1}})_{\mathbb{Z}_p}
\cong \mathbf{W}(n;\underline{1})_{\mathbb{Z}_p}
=\text{Span}_{\mathbb{Z}_p}\{x^{(\alpha)}D_i\mid 0\le \alpha\le
\tau, 1\le i\le n\}$. Moreover, we have $\mathbf{W}(n;\underline{1})
=\mathcal{K}\otimes_{\mathbb{Z}_p}\mathbf{W}(n;\underline{1})_{\mathbb{Z}_p}
=\mathcal{K}\mathbf{W}(n;\underline{1})_{\mathbb{Z}_p}$, and
$\mathbf{W}^+_{\mathcal{K}}=\mathcal{K}\mathbf{W}^+_{\mathbb{Z}_p}$.

Observe that the ideal
$J_{\underline{1}}:=\mathcal{K}(J_{\underline{1}})_{\mathbb{Z}_p}$
generates an ideal of $U(\mathbf{W}^+_{\mathcal{K}})$ over
$\mathcal{K}$, denoted by
$J:=J_{\underline{1}}U(\mathbf{W}^+_{\mathcal{K}})$, where
$\mathbf{W}^+_{\mathcal{K}}/J_{\underline{1}}\cong
\mathbf{W}(n;\underline{1})$. Based on the formulae (24) \& (25), we
see that $J$ is a Hopf ideal of $U(\mathbf{W}^+_{\mathcal{K}})$ such
that $U(\mathbf{W}^+_{\mathcal{K}})/J\cong
U(\mathbf{W}(n;\underline{1}))$. Note that the elements
$\frac{1}{\alpha!}x^\alpha D_i$ in $\mathbf{W}^+_{\mathcal{K}}$ for
$0\le\alpha\le\tau$ will be identified with $x^{(\alpha)}D_i$ in
$\mathbf{W}(n;\underline{1})$ and those in $J_{\underline{1}}$ with
$0$. Hence, by Corollary 2.10, we get the quantization of
$U(\mathbf{W}(n;\underline{1}))$ over
$U(\mathbf{W}(n;\underline{1}))[[t]]$ as follows.
\begin{theorem}
Given two distinguished elements $h:=x^{(\epsilon_k)}D_k$,
$e:=2x^{(2\epsilon_k)}D_k$ $(1\le k\le n)$ with $[h,\,e]=e$, the
corresponding quantization of $U(\mathbf{W}(n;\underline{1}))$ over
$U(\mathbf{W}(n;\underline{1}))[[t]]$ with the algebra structure
undeformed is given by
\begin{gather*}
\Delta({x^{(\alpha)}D_i}){=}x^{(\alpha)}D_i{\otimes}
(1{-}et)^{\alpha_k{-}\delta_{ik}}{+}\sum\limits_{\ell=0}^{p{-}1}{({-}1)}^\ell\bar
C_\ell\, h^{\lg
\ell\rg}{\otimes}(1{-}et)^{{-}\ell}x^{(\alpha{+}\ell\epsilon_k)}D_it^\ell,
\tag{33}\\
\end{gather*}
\begin{gather*}
S(x^{(\alpha)}D_i){=}{-}(1{-}et)^{-\alpha_k{+}\delta_{ik}}\cdot\Bigl(\sum\limits_{\ell=0}^{p{-}1}
\bar C_\ell\,x^{(\alpha{+}\ell\epsilon_k)}D_i\cdot h_1^{\lg
\ell\rg}t^\ell\Bigr),
\tag{34}\\
\varepsilon(x^{(\alpha)}D_i)=0, \tag{35}
\end{gather*}
where $0\le \alpha \le \tau$, $\bar C_\ell=\bar A_\ell-\bar B_\ell$,
$\bar A_\ell=\ell!\binom{\alpha_k{+}\ell}{\ell}A_\ell\,(\text{\rm
mod} \,p)$, $\bar B_\ell=\ell!
\binom{\alpha_k{+}\ell}{\ell}B_\ell\,(\text{\rm mod}\,p)$.
\end{theorem}

\subsection{Modulo ``$p$-restrictedness" reduction}
Assume that $I$ is the ideal of $U(\mathbf{W}(n;\underline{1}))$
over $\mathcal{K}$ generated by
$(x^{(\epsilon_i)}D_i)^p-x^{(\epsilon_i)}D_i$ and
$(x^{(\alpha)}D_i)^p$ with $\alpha\ne \epsilon_i$ for $0\le\alpha\le
\tau$ and $1\le i\le n$.
$\mathbf{u}(\mathbf{W}(n;\underline{1}))=U(\mathbf{W}(n;\underline{1}))/I$
is of dimension $p^{np^n}$. In order to get a reasonable
quantization of finite dimension for
$\mathbf{u}(\mathbf{W}(n;\underline{1}))$ in characteristic $p$, at
first, it is necessary to clarify (in concept) what is the
underlying vector space in which the required $t$-deformed object
exists. According to our modulo $p$ reduction approach, it should be
induced from the topologically free $\mathcal{K}[[t]]$-algebra
$U(\mathbf{W}(n;\underline{1}))[[t]]$ given in Theorem 3.1, or more
naturally, related with the topologically free
$\mathcal{K}[[t]]$-algebra
$\mathbf{u}(\mathbf{W}(n;\underline{1}))[[t]]$. As
$\mathbf{u}(\mathbf{W}(n;\underline{1}))$ is finite dimensional, as
algebras, we have (see p. 2, \cite{ES})
$$
\mathbf{u}(\mathbf{W}(n;\underline{1}))[[t]]\cong
\mathbf{u}(\mathbf{W}(n;\underline{1}))\otimes_{\mathcal{K}}
\mathcal{K}[[t]].\leqno(36)
$$
So the standard Hopf algebra structure
$(\mathbf{u}(\mathbf{W}(n;\underline{1}))[[t]],
m,\iota,\Delta_0,S_0,\varepsilon_0)$ can be viewed as the tensor
Hopf algebra structure of the standard ones on
$\mathbf{u}(\mathbf{W}(n;\underline{1}))$ and on $\mathcal{K}[[t]]$,
while the expected {\it twisted} Hopf algebra structure on
$(\mathbf{u}(\mathbf{W}(n;\underline{1}))[[t]], m,\iota$,
$\Delta,S,\varepsilon)$ induced from Theorem 3.1 when restricted to
the sub-Hopf algebra $\mathcal{K}[[t]]$ should still go back to the
original one on $\mathcal{K}[[t]]$ itself.

Keep the above viewpoint in mind. Furthermore, we need another
observation below.
\begin{lemm} $(\text{\rm i})$ \ $(1-et)^p\equiv 1 \quad (\text{\rm mod}\,p, I)$.

$(\text{\rm ii})$ \ $(1-et)^{-1}\equiv 1+et+\cdots+e^{p-1}t^{p-1}
\quad (\text{\rm mod}\,
 p, I)$.

$(\text{\rm iii})$ \ $h_a^{\lg \ell\rg} \equiv 0 \quad (\text{\rm
mod} \, p, I)$ for $\ell \geq p$, and $a\in\mathbb{Z}_p$.
\end{lemm}
\begin{proof} (i), (ii) follow from $e^p=0$ in $\mathbf{u}(\mathbf{W}(n;\underline{1}))$.

(iii)  For $\ell\in\mathbb{Z}_+$, there is a unique decomposition
$\ell=\ell_0+\ell_1p$ with $0\le \ell_0<p$ and $\ell_1\ge 0$.
Using the formulae (1) \& (3), we have
\begin{equation*}
\begin{split}
h_a^{\langle \ell\rangle}& =h_a^{\langle \ell_0\rangle}\cdot
h_{a+\ell_0}^{\langle \ell_1p\rangle}\equiv h_a^{\langle
\ell_0\rangle}\cdot (h_{a+\ell_0}^{\langle
p\rangle})^{\ell_1}\,\qquad
(\text{mod } p)\\
&\equiv h_a^{\langle \ell_0\rangle}\cdot (h^p-h)^{\ell_1}\quad
(\text{mod } p),
\end{split}
\end{equation*}
where we used the facts that $(x+1)(x+2)\cdots(x+p-1)\equiv
x^{p-1}-1\; (\text{mod } p)$, and $(x+a+\ell_0)^p\equiv
x^p+a+\ell_0\ (\text{mod } p)$. Hence, $h_a^{\lg \ell\rg} \equiv
0$ (mod $p, \,I$) for $\ell \geq p$.
\end{proof}
This Lemma, combining with Theorem 3.1, indicates that the required
$t$-deformation of $\mathbf{u}(\mathbf{W}(n;\underline{1}))$ (if it
exists) only happens in a $p$-truncated polynomial ring (with
degrees of $t$ less than $p$) with coefficients in
$\mathbf{u}(\mathbf{W}(n;\underline{1}))$. Let $\mathcal{K}[t]_p$ be
a $p$-truncated polynomial ring (of small possible dimension) over
$\mathcal{K}$, which should inherit a standard Hopf algebra
structure from that on $\mathcal{K}[[t]]$ with respect to modulo a
Hopf ideal of it. In the modular case, such a Hopf ideal in
$\mathcal{K}[[t]]$ has to take the form $(t^p-qt)$ generated by a
$p$-polynomial $t^p-qt$ of degree $p$ for a parameter
$q\in\mathcal{K}$. Denote by $\mathcal{K}[t]_p^{(q)}$ the
corresponding quotient ring. That is to say
$$
\mathcal{K}[t]_p^{(q)}\cong \mathcal{K}[[t]]/(t^p-qt),
\qquad\text{\it for }\ q\in\mathcal{K}.\leqno(37)
$$
Thereby, we obtain the underlying ring for our required
$t$-deformation of $\mathbf{u}(\mathbf{W}(n;\underline{1}))$ over
$\mathcal{K}[t]_p^{(q)}$, denoted by
$\mathbf{u}_{t,q}(\mathbf{W}(n;\underline{1}))$. Moreover, as
standard Hopf algebras,
$$\mathbf{u}_{t,q}(\mathbf{W}(n;\underline{1}))\cong
\mathbf{u}(\mathbf{W}(n;\underline{1}))\otimes_{\mathcal{K}}\mathcal{K}[t]_p^{(q)}.\leqno(38)$$
Hence,
$\dim_{\mathcal{K}}\mathbf{u}_{t,q}(\mathbf{W}(n;\underline{1}))
=p\cdot\dim_{\mathcal{K}}\mathbf{u}(\mathbf{W}(n;\underline{1}))
=p^{1+np^n}$.

\begin{defi} With notations as above. A Hopf algebra
$(\mathbf{u}_{t,q}(\mathbf{W}(n;\underline{1}))$, $m,
\iota,\Delta,S,\varepsilon)$ is said to be a quantization of
$\mathbf{u}(\mathbf{W}(n;\underline{1}))$ (in characteristic $p$) if
it is a twisting of the standard Hopf algebra structure (as a tensor
Hopf algebra of the standard ones on
$\mathbf{u}(\mathbf{W}(n;\underline{1}))$ and on
$\mathcal{K}[t]_p^{(q)}$) with
$\mathbf{u}_{t,q}(\mathbf{W}(n;\underline{1}))/t\mathbf{u}_{t,q}(\mathbf{W}(n;\underline{1}))
$ $\cong \mathbf{u}(\mathbf{W}(n;\underline{1}))$.
\end{defi}
To describe $\mathbf{u}_{t,q}(\mathbf{W}(n;\underline{1}))$
explicitly, we still need an auxiliary Lemma.
\begin{lemm} $(\text{\rm i})$ \
$d^{(\ell)}(x^{(\alpha)}D_i)=\bar C_\ell
x^{(\alpha{+}\ell\epsilon_k)}D_i$, where
$d^{(\ell)}=\frac{1}{\ell!}(\text{\rm ad}\,e)^\ell$,
$e=2x^{(2\epsilon_k)}D_k$, and $\bar C_\ell$ as in Theorem 3.1.

$(\text{\rm ii})$ \
$d^{(\ell)}(x^{(\epsilon_i)}D_i)=\delta_{\ell,0}x^{(\epsilon_i)}D_i-\delta_{1,\ell}\delta_{ik}e$.

$(\text{\rm iii})$ \
$d^{(\ell)}((x^{(\alpha)}D_i)^p)=\delta_{\ell,0}(x^{(\alpha)}D_i)^p-\delta_{1,\ell}
\delta_{ik}\delta_{\alpha,\epsilon_i}e$.
\end{lemm}
\begin{proof} (i) For $0\le \alpha\le\tau$, by (19) and Theorem 3.1
(noting $B_\ell=\delta_{ik}A_{\ell{-}1}$), we get
\begin{equation*}
\begin{split}
d^{(\ell)}(x^{(\alpha)}D_i)&=\frac{1}{\alpha!}d^{(\ell)}(x^{\alpha-\epsilon_i}\partial_i)
=\frac{1}{\alpha!}x^{\alpha-\epsilon_i+\ell\epsilon_k}(A_\ell\partial_i-B_\ell\partial_k)\\
&=\ell!\dbinom{\alpha_k{+}\ell}{\ell}C_\ell\,x^{(\alpha{+}\ell\epsilon_k)}D_i\\
&=\bar C_\ell\,x^{(\alpha{+}\ell\epsilon_k)}D_i.
\end{split}
\end{equation*}

(ii) Noting that
$$\bar A_\ell=\ell! \binom{\alpha_k{+}\ell}{\ell}A_\ell
=\binom{\alpha_k{+}\ell}{\ell}\prod_{j=0}^{\ell{-}1}(\alpha_k{-}\delta_{ik}{+}j),\qquad
\bar B_\ell=\ell!\binom{\alpha_k{+}\ell}{\ell}\delta_{ik}
A_{\ell{-}1},$$ we obtain the following

If $i=k$: $\bar A_0=1$, $\bar A_\ell=0$ for $\ell\ge1$ and, $\bar
B_0=0$, $\bar B_1=2$, $\bar B_\ell=0$ for $\ell\ge 2$. So by Lemma
3.4 (i), we have
$d^{(1)}(x^{(\epsilon_k)}D_k)=-2x^{(2\epsilon_k)}D_k=-e$, and
$d^{(\ell)}(x^{(\epsilon_k)}D_k)=0$ for $\ell\ge2$.

If $i\ne k$: $\bar A_0=1$, $\bar A_\ell=0$ for $\ell\ge 1$, and
$\bar B_\ell=0$ for $\ell\ge0$, namely,
$d^{(\ell)}(x^{(\epsilon_i)}D_i)=0$ for $\ell\ge1$.

In both cases, we arrive at the result as required.

\smallskip
(iii) From (12), we obtain that for $0\le \alpha\le\tau$,
\begin{equation*}
\begin{split}
d^{(1)}\,((x^{(\alpha)}D_i)^p)&=\frac{1}{(\alpha!)^p}[\,e,(x^{\alpha}D_i)^p\,]=\frac{1}{\alpha!}[\,e,(x^{\alpha-\epsilon_i}\partial_i)^p\,]\\
&=\frac{1}{\alpha!}\sum\limits_{\ell=1}^p(-1)^\ell\dbinom{p}
{\ell}(x^{\alpha-\epsilon_i}\partial_i)^{p-\ell}\cdot
x^{\epsilon_k{+}\ell(\alpha{-}\epsilon_i)}(a_\ell
\partial_k-b_\ell\partial_i)\\
&\equiv
-\frac{a_p}{\alpha!}\,x^{2\epsilon_k{+}p(\alpha{-}\epsilon_i)}
D_k\qquad(\text{mod }\,p\,)\\
&\equiv \begin{cases} -{a_p}\,e,\qquad & \text{\it if }\quad \alpha=\epsilon_i\\
0,\qquad & \text{\it if }\quad \alpha\ne\epsilon_i
\end{cases}\qquad(\text{mod }\,J),
\end{split}
\end{equation*}
where we get the last ``$\equiv$" by using the identification with
respect to modulo the ideal $J$ as before, and
$a_\ell=\prod\limits_{j=0}^{\ell-1}(\delta_{ik}{+}j(\alpha_i{-}1)),\
b_\ell=\ell\,(\alpha_k-\delta_{ik})a_{\ell-1}$, and
$a_p=\delta_{ik}$
for $\alpha=\epsilon_i$. 

Consequently, by definition of $d^{(\ell)}$, we obtain
$d^{(\ell)}((x^{(\alpha)}D_i)^p)=0$ in
$\mathbf{u}(\mathbf{W}(n;\underline{1}))$ for $2\le \ell\le p-1$ and
any $0\le\alpha\le\tau$.
\end{proof}

Based on Theorem 3.1, Definition 3.3 and Lemma 3.4, we arrive at
\begin{theorem} With the specific choice of the two distinguished elements
$h:=x^{(\epsilon_k)}D_k$, $e:=2x^{(2\epsilon_k)}D_k$ $(1\le k\le n)$
with $[h,\,e]=e$, there is a noncommutative and noncocummtative Hopf
algebra
$(\mathbf{u}_{t,q}(\mathbf{W}(n;\underline{1})),m,\iota,\Delta,S,\varepsilon)$
over $\mathcal{K}[t]_p^{(q)}$ with its algebra structure undeformed,
whose coalgebra structure is given by
\begin{gather*}
\Delta(x^{(\alpha)}D_i)=x^{(\alpha)}D_i\otimes
(1-et)^{\alpha_k{-}\delta_{ik}}
+\sum\limits_{\ell=0}^{p-1}(-1)^\ell h^{\lg \ell\rg}\otimes
(1-et)^{-\ell}\cdot d^{(\ell)}(x^{(\alpha)}D_i)t^\ell,
\tag{39}\\
S(x^{(\alpha)}D_i)=-(1-et)^{\delta_{ik}{-}\alpha_k}\Bigl(\sum\limits_{\ell=0}^{p-1}d^{(\ell)}(x^{(\alpha)}D_i)\cdot
h_1^{\lg
\ell\rg}t^\ell\Bigr),\tag{40}\\
\varepsilon(x^{(\alpha)}D_i)=0,  \tag{41}
\end{gather*}
where $0\le\alpha\le\tau$. It is finite dimensional and
$\dim_{\mathcal{K}}\mathbf{u}_{t,q}(\mathbf{W}(n;\underline{1}))=p^{1{+}np^n}$.
\end{theorem}
\begin{proof} Let $I_t$ be the ideal of $(U(\mathbf{W}(n;\underline{1}))[[t]],
m, \iota,\Delta,S,\varepsilon)$ generated by $I$ and $t^p-qt$
($q\in\mathcal{K}$). In what follows, we shall show that the ideal
$I_t$ is a Hopf ideal of the {\it twisted} Hopf algebra
$U(\mathbf{W}(n;\underline{1}))[[t]]$ given in Theorem 3.1. To this
end, it suffices to verify that $\Delta$ and $S$ preserve the
elements in $I$ since $\Delta(t^p-qt)=(t^p-qt)\otimes 1+1\otimes
(t^p-qt)$ and $S(t^p-qt)=-(t^p-qt)$.

\smallskip
(I) \ By Lemmas 2.8, 3.2 \& 3.4 (iii), we obtain
\begin{equation*}
\begin{split}
\Delta((x^{(\alpha)}D_i)^p) &=(x^{(\alpha)}D_i)^p\otimes
(1{-}et)^{p(\alpha_k{-}\delta_{ik})}\\
&\quad+\sum\limits_{\ell=0}^{\infty} ({-}1)^\ell h^{\lg
\ell\rg}\otimes
(1{-}et)^{{-}\ell}d^{(\ell)}((x^{(\alpha)}D_i)^p)t^\ell\
\\
&\equiv(x^{(\alpha)}D_i)^p{\otimes}1+\sum\limits_{\ell=0}^{p{-}1}
({-}1)^\ell h^{\lg \ell\rg}{\otimes}
(1{-}et)^{{-}\ell}d^{(\ell)}((x^{(\alpha)}D_i)^p)t^\ell\\
&\qquad\quad (\text{\rm mod }\, p,\; I_t\otimes
U(\mathbf{W}(n;\underline{1}))[[t]]+U(\mathbf{W}(n;\underline{1}))[[t]]\otimes
I_t)
\\
&=(x^{(\alpha)}D_i)^p{\otimes}1+1{\otimes}(x^{(\alpha)}D_i)^p
+h{\otimes}(1{-}et)^{-1}\delta_{ik}\delta_{\a,\epsilon_i}et.\\
&\qquad\qquad (\text{\rm mod }\, I_t\otimes
U(\mathbf{W}(n;\underline{1}))[[t]]+U(\mathbf{W}(n;\underline{1}))[[t]]\otimes
I_t)
\end{split}\tag{42}
\end{equation*}
Hence, when $\alpha\ne\epsilon_i$, we get
\begin{equation*}
\begin{split}
\Delta((x^{(\alpha)}D_i)^p)&\equiv(x^{(\alpha)}D_i)^p\otimes
1+1\otimes
(x^{(\alpha)}D_i)^p\\
&\subseteq I_t\otimes
U(\mathbf{W}(n;\underline{1}))[[t]]+U(\mathbf{W}(n;\underline{1}))[[t]]\otimes
I_t.
\end{split}
\end{equation*}
When $\al=\epsilon_i$, by Lemma 3.4 (ii), (27) becomes
$$\Delta(x^{(\epsilon_i)}D_i)=x^{(\epsilon_i)}D_i\otimes 1+
1\otimes x^{(\epsilon_i)}D_i+\delta_{ik}h\otimes (1-et)^{-1}et.$$
Combining with (42), we obtain
\begin{equation*}
\begin{split}
\Delta((x^{(\epsilon_i)}D_i)^p-x^{(\epsilon_i)}D_i)&\equiv((x^{(\epsilon_i)}D_i)^p-
x^{(\epsilon_i)}D_i)\otimes
1+1\otimes ((x^{(\epsilon_i)}D_i)^p-x^{(\epsilon_i)}D_i)\\
&\subseteq I_t\otimes
U(\mathbf{W}(n;\underline{1}))[[t]]+U(\mathbf{W}(n;\underline{1}))[[t]]\otimes
I_t.
\end{split}
\end{equation*}

Thereby, we prove that the ideal $I_t$ is also a coideal of the Hopf
algebra $ U(\mathbf{W}(n;\underline{1}))[[t]]$.

\smallskip
(II) \ By Lemmas 2.8, 3.2 \& 3.4 (iii), we have
\begin{equation*}
\begin{split}
S((x^{(\alpha)}D_i)^p) &=-(1{-}et)^{-p(\alpha_k-\delta_{ik})}
\sum\limits_{\ell=0}^{\infty} d^{(\ell)}((x^{(\alpha)}D_i)^p)\cdot
h_1^{\lg
\ell\rg}t^\ell\\
&\equiv -(x^{(\alpha)}D_i)^p-\sum\limits_{\ell=1}^{p-1}
d^{(\ell)}((x^{(\alpha)}D_i)^p)\cdot h_1^{\lg \ell\rg}t^\ell \quad
(\text{mod }(p, I))\\
&=-(x^{(\alpha)}D_i)^p+\delta_{ik}\delta_{\alpha,\epsilon_i}e\cdot
h_1^{\lg 1\rg} t.
\end{split}\tag{43}
\end{equation*}
Hence, when $\alpha\ne\epsilon_i$, we get
$$S\bigl((x^{(\al)}D_i)^p\bigr)\equiv-(x^{(\al)}D_i)^p\equiv 0\quad (\text{\rm mod} \, I_t).$$
When $\al=\epsilon_i$, by Lemma 3.4 (ii), (28) reads as
$S(x^{(\epsilon_i)}D_i)=-x^{(\epsilon_i)}D_i+\delta_{ik}e\cdot
h_1^{\langle 1\rangle}t$. Combining with (43), we obtain
$$S\bigl((x^{(\epsilon_i)}D_i)^p-x^{(\epsilon_i)}D_i\bigr)\equiv-\bigl((x^{(\epsilon_i)}D_i)^p-x^{(\epsilon_i)}D_i\bigr)\equiv
0\quad (\text{\rm mod} \, I_t).$$

Thereby, we show that the ideal $I_t$ is indeed preserved by the
antipode $S$ of the quantization
$U(\mathbf{W}(n;\underline{1}))[[t]]$ given in Theorem 3.1.

\smallskip
(III) It is obvious to notice that
$\varepsilon((x^{(\alpha)}D_i)^p)=0$ for all $0\le\alpha\le\tau$.

\smallskip
In other words, we prove that $I_t$ is a Hopf ideal in
 $U(\mathbf{W}(n;\underline{1}))[[t]]$. We thus obtain the required
 $t$-deformation on
 $\mathbf{u}_{t,q}(\mathbf{W}(n;\underline{1}))$ for the Cartan type
 simple modular restricted Lie algebra of $W$ type
 --- the Jacobson-Witt algebra
$\mathbf{W}(n;\underline{1})$.
\end{proof}

\begin{remark} (i) \
If we set $f=(1-et)^{-1}$, then by Lemma 3.4, Theorem 3.5, we have
$$[h,f]=f^2-f,\quad h^p=h, \quad f^p=1, \quad
\Delta(h)=h\otimes f+1\otimes h,$$ where $f$ is a group-like
element, and $S(h)=-hf^{-1}$, $\varepsilon(h)=0$. So the subalgebra
generated by $h$ and $f$ is a Hopf subalgebra of
$\mathbf{u}_{t,q}(\mathbf{W}(n;\underline{1}))$, which is isomorphic
to the well-known Radford Hopf algebra over $\mathcal{K}$ in
characteristic $p$.

\smallskip
(ii) \ According to our argument above, given a parameter
$q\in\mathcal{K}$, one can specialize $t$ to any one of roots of the
$p$-polynomial $t^p-qt\in\mathcal{K}[t]$ in a split field of
$\mathcal{K}$. For instance, if take $q=1$, then one can specialize
$t$ to any scalar in $\mathbb{Z}_p$. If set $t=0$, then we get the
original standard Hopf algebra structure of
$\mathbf{u}(\mathbf{W}(n;\underline{1}))$. In this way, we indeed
get a new Hopf algebra structure over the same restricted universal
enveloping algebra $\mathbf{u}(\mathbf{W}(n;\underline{1}))$ over
$\mathcal{K}$ under the assumption that $\mathcal{K}$ is
algebraically closed, which has the new coalgebra structure induced
by Theorem 3.5, but has dimension $p^{np^n}$.
\end{remark}

\subsection{More quantizations }
Carrying out the modular reduction process for those pairwise
different products of some basic Drinfel'd twists as stated in
Remark 2.11, we will get many more new families of noncommutative
and noncocommutative Hopf algebras of dimension $p^{1{+}np^n}$ with
indeterminate $t$ or of dimension $p^{np^n}$ with specializing $t$
into a scalar in $\mathcal{K}$. We have the following general
results about the quantizations under concern.

Maintain the notations as in Theorem 2.14 and set
$d_I^{(\eta\ell)}=\prod\limits_{k \in I}d_k^{( \eta_k\ell_k)}$,
where $d_k^{(\ell_k)}=\frac{1}{\ell_k!}(\text{\rm
ad}\,e(k))^{\ell_k}$. Then we have
\begin{theorem}
For each given $\eta=(\eta_1,\cdots,\eta_n)$ with
$\eta_i\in\{0,1\}$, there exists a noncommutative and
noncocummtative Hopf algebra
$(\mathbf{u}_{t,q}(\mathbf{W}(n;\underline{1})),m,\iota,$
$\Delta,S,\varepsilon)$ over $\mathcal{K}[t]_p^{(q)}$ with the
algebra structure undeformed, whose coalgebra structure is given by
\begin{gather*}
\Delta(x^{(\alpha)}D_i)=x^{(\alpha)}D_i{\otimes}
(1{-}e(I)t)^{\eta(\alpha{-}\delta_{iI})}
{+}\sum\limits_{\eta\ell=\underline{0}}^{\eta\tau}{({-}1)}^{|\eta\ell|}
h(I)^{\lg \eta\ell\rg}{\otimes}\tag{44}
\\
(1{-}e(I)t)^{{-}\eta\ell}\cdot d_I^{(\eta
\ell)}(x^{(\alpha)}D_i)t^{|\eta\ell|},
\\
S(x^{(\alpha)}D_i)={-}(1{-}e(I)t)^{-\eta(\alpha{-}\delta_{iI})}
\Bigl(\sum\limits_{\eta\ell=\underline{0}}^{\eta\tau} d_I^{(\eta
\ell)} (x^{(\alpha)}D_i )h(I)_1^{\lg
\eta\ell\rg}t^{|\eta\ell|}\Bigr), \tag{45} \\
\varepsilon(x^{(\alpha)}D_i)=0,  \tag{46}
\end{gather*}
where $0\leq \alpha \leq \tau$. It is finite dimensional and
$\dim_{\mathcal{K}}\mathbf{u}_{t,q}(\mathbf{W}(n;\underline{1}))=p^{1{+}np^n}$.
\end{theorem}
For the proof, it suffices to work with the case in Lemma 2.13. To
this end, we begin to perform the modular reduction process for the
quantization integral form corresponding to the Drinfel'd twist
$\mathcal{F}=\mathcal{F}(m)\mathcal{F}(k)$ ($m\ne k$).

\begin{lemm}
Given two pairs of distinguished elements
$h(k):=x^{(\epsilon_k)}D_k$, $e(k):=2x^{(2\epsilon_k)}D_k$ with
$[h(k),\,e(k)]=e(k)$ and $h(m):=x^{(\epsilon_m)}D_m$,
$e(m):=2x^{(2\epsilon_m)}D_m$ $(1\le m\neq k\le n)$ with
$[h(m),\,e(m)]=e(m)$, the corresponding quantization of
$U(\mathbf{W}(n;\underline{1}))$ over
$U(\mathbf{W}(n;\underline{1}))[[t]]$  with the algebra structure
undeformed is given by
\begin{gather*}
\Delta(x^{(\alpha)}D_i)=x^{(\alpha)}D_i{\otimes}
(1{-}e(k)t)^{\alpha_k{-}\delta_{ik}}(1{-}e(m)t)^{\alpha_m{-}\delta_{im}}
{+}\sum\limits_{\ell,n=0}^{p-1}{({-}1)}^{\ell+n}\tag{47} \\
\cdot \bar C(k)_\ell \bar C(m)_n h(k)^{\lg \ell\rg} h(m)^{\lg n \rg}
{\otimes} (1{-}e(k)t)^{{-}\ell} (1{-}e(m)t)^{{-}n}\cdot
x^{(\alpha{+}\ell\epsilon_k+n\epsilon_m)}D_it^{\ell+n},\\
S(x^{(\alpha)}D_i)={-}(1{-}e(k)t)^{-\alpha_k{+}\delta_{ik}}(1{-}e(m)t)^{-\alpha_m{+}\delta_{im}}
 \Bigl(\sum\limits_{\ell,n=0}^{p-1}\bar C(k)_\ell\bar C(m)_n \tag{48} \\ \cdot
x^{(\alpha{+}\ell\epsilon_k+n\epsilon_m)}D_i\cdot h(m)_1^{\lg
n\rg}h(k)_1^{\lg \ell\rg}t^{\ell+n}\Bigr),\\
\varepsilon(x^{(\alpha)}D_i)=0,  \tag{49}
\end{gather*}
where  $0\le \alpha \le \tau$, $\bar C(k)_\ell=\bar A(k)_\ell-\bar
B(k)_\ell$, $\bar
A(k)_\ell=\ell!\binom{\alpha_k{+}\ell}{\ell}A(k)_\ell\,(\text{\rm
mod} \,p)$, $\bar
B(k)_\ell=\ell!\binom{\alpha_k{+}\ell}{\ell}B(k)_\ell\,(\text{\rm
mod}\,p)$.
\end{lemm}
We have the following two Lemmas about the quantization of
$U(\mathbf{W}(n;\underline{1}))$ over
$U(\mathbf{W}(n;\underline{1}))[[t]]$ defined in Lemma 3.8.
\begin{lemm} For $s\ge 1$, one has
\begin{gather*}
\Delta((x^\alpha D_i)^s)=\sum_{0\le j\le s\atop n, \ell\ge
0}\dbinom{s}{j}({-}1)^{n+\ell}(x^{\alpha}D_i)^jh(k)^{\langle
\ell\rangle}h(m)^{\langle n\rangle}\tag{\text{\rm i}}\otimes\\
(1{-}e(k)t)^{j(\al_k-\delta_{ik}){-}\ell}(1{-}e(m)t)^{j(\al_m-\delta_{im}){-}n}
d_m^{(n)}d_k^{(\ell)}((x^\al D_i)^{s{-}j})t^{\ell+n}.\\
S((x^{\alpha} D_i)^s)=
(-1)^s(1-e(k)t)^{-s(\al_k-\delta_{ik})}(1-e(m)t)^{-s(\al_m-\delta_{im})}
\tag{\text{\rm ii}} \\ \cdot\Bigl(\sum\limits_{n,\ell=0}^{\infty}
d_m^{(n)}d_k^{(\ell)}((x^{\alpha}D_i)^s)\cdot h(k)_1^{\lg
\ell\rg}h(m)_1^{\lg n\rg}t^{n+\ell}\Bigr).
\end{gather*}
\end{lemm}
\begin{proof} By Lemma 2.8, (18), (14) and Lemma 2.2, we obtain
\begin{equation*}
\begin{split}
\Delta((x^{\alpha}D_i)^s)&=\mathcal{F}\Bigl(x^{\al}D_i\otimes
1+1\otimes x^{\alpha}D_i\Bigr)^s\mathcal{F}^{-1}\\
&=\mathcal{F}(m)\Bigl(\sum_{0\le j\le s\atop
\ell\ge0}\dbinom{s}{j}({-}1)^\ell(x^{\alpha}D_i)^jh(k)^{\langle
\ell\rangle}\otimes(1{-}e(k)t)^{j(\al_k-\delta_{ik}){-}\ell}\\
& \quad
\cdot d_k^{(\ell)}((x^\al D_i)^{s{-}j})t^\ell\Bigr)\mathcal{F}(m)^{-1}\\
&=\mathcal{F}(m)\Bigl(\sum_{0\le j\le s\atop
\ell\ge0}\dbinom{s}{j}({-}1)^\ell((x^{\alpha}D_i)^j\otimes 1)
(h(k)^{\langle
\ell\rangle}\otimes(1{-}e(k)t)^{j(\al_k-\delta_{ik}){-}\ell})\\
& \quad
\cdot (1\otimes d_k^{(\ell)}((x^\al D_i)^{s{-}j})t^\ell)\Bigr)\mathcal{F}(m)^{-1}\\
&=\mathcal{F}(m)\sum_{0\le j\le s\atop n, \ell\ge
0}\dbinom{s}{j}({-}1)^{n+\ell}((x^{\alpha}D_i)^j\otimes
1)\mathcal{F}(m)_n^{-1}h(k)^{\langle \ell\rangle}h(m)^{\langle
n\rangle}\\  & \quad
\otimes(1{-}e(k)t)^{j(\al_k-\delta_{ik}){-}\ell}
d_m^{(n)}d_k^{(\ell)}((x^\al D_i)^{s{-}j})t^{\ell+n}\\
&=\sum_{0\le j\le s\atop n, \ell\ge
0}\dbinom{s}{j}({-}1)^{n+\ell}\mathcal{F}(m)\mathcal{F}(m)_{n-j(\al_m-\delta_{im})}^{-1}
((x^{\alpha}D_i)^j\otimes 1)h(k)^{\langle \ell\rangle}h(m)^{\langle
n\rangle}\\  & \quad
\otimes(1{-}e(k)t)^{j(\al_k-\delta_{ik}){-}\ell}
d_m^{(n)}d_k^{(\ell)}((x^\al D_i)^{s{-}j})t^{\ell+n}\\
&=\sum_{0\le j\le s\atop n, \ell\ge
0}\dbinom{s}{j}({-}1)^{n+\ell}(x^{\alpha}D_i)^jh(k)^{\langle
\ell\rangle}h(m)^{\langle
n\rangle}\otimes(1{-}e(k)t)^{j(\al_k-\delta_{ik}){-}\ell}
\\  & \quad \cdot(1{-}e(m)t)^{j(\al_m-\delta_{im}){-}n}
d_m^{(n)}d_k^{(\ell)}((x^\al D_i)^{s{-}j})t^{\ell+n}.
\end{split}
\end{equation*}
Again by (17) and Lemma 2.2,
\begin{equation*}
\begin{split}
S((x^{\alpha}D_i)^s)&=u^{-1}S_0((x^{\alpha}\partial)^s)\,u=(-1)^s
v\cdot (x^{\alpha}\partial)^s\cdot u \\ &=(-1)^s
v(m)\Bigl((1-e(k)t)^{-s(\al_k-\delta_{ik})}\cdot\Bigl(\sum\limits_{\ell=0}^{\infty}
d_k^{(\ell)}((x^{\alpha}D_i)^s)\cdot h(k)_1^{\lg
\ell\rg}t^\ell\Bigr)\Bigr)u(m)
\\ &=(-1)^s
v(m)u(m)_{s(\al_m-\delta_{im})}(1-e(k)t)^{-s(\al_k-\delta_{ik})}\\
&\quad \cdot \Bigl(\sum\limits_{n,\ell=0}^{\infty}
d_m^{(n)}d_k^{(\ell)}((x^{\alpha}D_i)^s)\cdot h(k)_1^{\lg \ell\rg}
h(m)_1^{\lg n \rg}t^{n+\ell}\Bigr) \\
\end{split}
\end{equation*}
\begin{equation*}
\begin{split}
&=(-1)^s
(1-e(m)t)^{-s(\al_m-\delta_{im})}(1-e(k)t)^{-s(\al_k-\delta_{ik})}\\
& \quad \cdot \Bigl(\sum\limits_{n,\ell=0}^{\infty}
d_m^{(n)}d_k^{(\ell)}((x^{\alpha}D_i)^s)\cdot h(k)_1^{\lg \ell\rg}
h(m)_1^{\lg n \rg}t^{n+\ell}\Bigr).
\end{split}
\end{equation*}
So, the proof is complete.
\end{proof}

Lemma 3.4 implies the following
\begin{lemm} $(\text{\rm i})$ \
$d_m^{(n)}d_k^{(\ell)}(x^{(\alpha)}D_i)=\bar C(k)_\ell\bar C(m)_n
x^{(\alpha{+}\ell\epsilon_k+n\epsilon_m)}D_i$, where
$d_k^{(\ell)}=\frac{1}{\ell!}(\text{\rm ad}\,e(k))^\ell$,
$e(k)=2x^{(2\epsilon_k)}D_k$, and $\bar C(k)_\ell$, $\bar C(m)_n$ as
in Lemma 3.8.

$(\text{\rm ii})$ \
$d_m^{(n)}d_k^{(\ell)}(x^{(\epsilon_i)}D_i)=\delta_{\ell,0}\delta_{n,0}x^{(\epsilon_i)}D_i
-(\delta_{n,0}\delta_{1,\ell}\delta_{ik}+\delta_{\ell,0}\delta_{1,n}\delta_{im})\,e(i)
$.

$(\text{\rm iii})$ \
$d_m^{(n)}d_k^{(\ell)}((x^{(\alpha)}D_i)^p)=\delta_{\ell,0}\delta_{n,0}(x^{(\alpha)}D_i)^p-(\delta_{n,0}\delta_{1,\ell}
\delta_{ik}+\delta_{\ell,0}\delta_{1,n}
\delta_{im})\,\delta_{\alpha,\epsilon_i}e(i)$.
\end{lemm}

Using Lemmas 3.2, 3.4, 3.9 \& 3.10, we get a new Hopf algebra
structure over the same algebra
$\mathbf{u}_{t,q}(\mathbf{W}(n;\underline{1}))$ as follows.
\begin{theorem} Given two pairs of distinguished elements
$h(k):=x^{(\epsilon_k)}D_k$, $e(k):=2x^{(2\epsilon_k)}D_k$ with
$[h(k),\,e(k)]=e(k);$ $h(m):=x^{(\epsilon_m)}D_m$,
$e(m):=2x^{(2\epsilon_m)}D_m$ $(1\le m\neq k\le n)$ with
$[h(m),\,e(m)]=e(m)$, there is a noncommutative and noncocummtative
Hopf algebra
$(\mathbf{u}_{t,q}(\mathbf{W}(n;\underline{1})),m,\iota,$
$\Delta,S,\varepsilon)$ over $\mathcal{K}[t]_p^{(q)}$ with its
algebra structure undeformed, whose coalgebra structure is given by
\begin{gather*}
\Delta(x^{(\alpha)}D_i)=x^{(\alpha)}D_i\otimes
(1-e(k)t)^{\alpha_k{-}\delta_{ik}}(1-e(m)t)^{\alpha_m{-}\delta_{im}}
+\sum\limits_{n,\ell=0}^{p-1}(-1)^{\ell+n} \tag{50}\\ \cdot
h(k)^{\lg \ell\rg} h(m)^{\lg n \rg}\otimes
(1-e(k)t)^{-\ell}(1-e(m)t)^{-n}
d_k^{(\ell)}d_m^{(n)}(x^{(\alpha)}D_i)t^{\ell+n},\\
S(x^{(\alpha)}D_i)=-(1-e(k)t)^{\delta_{ik}{-}\alpha_k}(1-e(m)t)^{\delta_{im}{-}\alpha_m} \tag{51}\\
\cdot
\Bigl(\sum\limits_{n,\ell=0}^{p-1}d_k^{(\ell)}d_m^{(n)}(x^{(\alpha)}D_i)
 h(k)_1^{\lg \ell\rg}h(m)_1^{\lg n\rg}t^{\ell+n}\Bigr),\\
\varepsilon(x^{(\alpha)}D_i)=0,  \tag{52}
\end{gather*}
where $0\le\alpha\le\tau$. It is finite dimensional and
$\dim_{\mathcal{K}}\mathbf{u}_{t,q}(\mathbf{W}(n;\underline{1}))=p^{1{+}np^n}$.
\end{theorem}
\begin{proof} Let $I_t$ be the ideal of $(U(\mathbf{W}(n;\underline{1}))[[t]],
m, \iota,\Delta,S,\varepsilon)$ generated by $I$ and $t^p-qt$
($q\in\mathcal{K}$). We begin to show that the ideal $I_t$ is a Hopf
ideal of the {\it twisted} Hopf algebra
$U(\mathbf{W}(n;\underline{1}))[[t]]$ given in Lemma 3.8. To this
purpose, it suffices to verify that $\Delta$ and $S$ preserve the
elements in $I$ since $\Delta(t^p-qt)=(t^p-qt)\otimes 1+1\otimes
(t^p-qt)$ and $S(t^p-qt)=-(t^p-qt)$.

\smallskip
(I) \ By Lemmas 3.9, 3.2, 3.4 \& 3.10, we obtain
\begin{equation*}
\begin{split}
\Delta((x^{(\alpha)}D_i)^p) &{=}(x^{(\alpha)}D_i)^p\otimes
(1{-}e(k)t)^{p(\alpha_k{-}\delta_{ik})}(1{-}e(m)t)^{p(\alpha_m{-}\delta_{im})}\\
&\quad+\sum\limits_{n,\ell=0}^{\infty} ({-}1)^{n+\ell} h(k)^{\lg
\ell\rg} h(m)^{\lg n \rg}\otimes
(1{-}e(k)t)^{{-}\ell}(1{-}e(m)t)^{{-}n}\\
&\quad \cdot d_m^{(n)}d_k^{(\ell)}((x^{(\alpha)}D_i)^p)t^{n+\ell}
\\
\end{split}\tag{53}
\end{equation*}
\begin{equation*}
\begin{split}&{\equiv}(x^{(\alpha)}D_i)^p{\otimes}1+\sum\limits_{n,\ell=0}^{p{-}1}
({-}1)^{n+\ell} h(k)^{\lg \ell\rg} h(m)^{\lg n \rg}{\otimes}
(1{-}e(k)t)^{{-}\ell}(1{-}e(m)t)^{{-}n}\\
&\quad \cdot d_m^{(n)}d_k^{(\ell)}((x^{(\alpha)}D_i)^p)t^{n{+}\ell}\\
&\qquad\qquad (\text{\rm mod }\, p,\; I_t\otimes
U(\mathbf{W}(n;\underline{1}))[[t]]+U(\mathbf{W}(n;\underline{1}))[[t]]\otimes
I_t)\\
&{=}(x^{(\alpha)}D_i)^p{\otimes}1{+}\sum\limits_{n,\ell=0}^{p{-}1}
({-}1)^{n{+}\ell}h(k)^{\lg \ell\rg} h(m)^{\lg n \rg}{\otimes}
(1{-}e(k)t)^{{-}\ell} (1{-}e(m)t)^{{-}n}\\
& \quad \cdot\bigl(\delta_{\ell,0}\delta_{n,0}(x^{(\alpha)}D_i)^p
-(\delta_{n,0}\delta_{1,\ell}
\delta_{ik}+\delta_{\ell,0}\delta_{1,n}
\delta_{im})\,\delta_{\alpha,\epsilon_i}\,e(i)\bigr)t^{n+\ell}
\\
&\qquad\qquad (\text{\rm mod }\, I_t\otimes
U(\mathbf{W}(n;\underline{1}))[[t]]+U(\mathbf{W}(n;\underline{1}))[[t]]\otimes
I_t) \\
&{=}(x^{(\alpha)}D_i)^p{\otimes}1+1{\otimes}(x^{(\alpha)}D_i)^p
+h(k){\otimes}(1{-}e(k)t)^{-1}\delta_{ik}\delta_{\a,\epsilon_i}e(i)t\\
&\quad
+h(m){\otimes}(1{-}e(m)t)^{-1}\delta_{im}\delta_{\a,\epsilon_i}e(i)t.\\
&\qquad\qquad (\text{\rm mod }\, I_t\otimes
U(\mathbf{W}(n;\underline{1}))[[t]]+U(\mathbf{W}(n;\underline{1}))[[t]]\otimes
I_t).
\end{split}
\end{equation*}
Hence, when $\alpha\ne\epsilon_i$, we get
\begin{equation*}
\begin{split}
\Delta((x^{(\alpha)}D_i)^p)&\equiv(x^{(\alpha)}D_i)^p\otimes
1+1\otimes
(x^{(\alpha)}D_i)^p\\
&\subseteq I_t\otimes
U(\mathbf{W}(n;\underline{1}))[[t]]+U(\mathbf{W}(n;\underline{1}))[[t]]\otimes
I_t.
\end{split}
\end{equation*}
When $\al=\epsilon_i$, by Lemmas 3.4 and 3.10, (47) becomes
$\Delta(x^{(\epsilon_i)}D_i)=x^{(\epsilon_i)}D_i\otimes 1+ 1\otimes
x^{(\epsilon_i)}D_i+\delta_{ik}h(k)\otimes
(1-e(k)t)^{-1}e(i)t+\delta_{im}h(m)\otimes (1-e(m)t)^{-1}e(i)t.$
Combining with (53), we obtain
\begin{equation*}
\begin{split}
\Delta((x^{(\epsilon_i)}D_i)^p-x^{(\epsilon_i)}D_i)&\equiv((x^{(\epsilon_i)}D_i)^p
-x^{(\epsilon_i)}D_i)\otimes
1+1\otimes ((x^{(\epsilon_i)}D_i)^p-x^{(\epsilon_i)}D_i)\\
&\subseteq I_t\otimes
U(\mathbf{W}(n;\underline{1}))[[t]]+U(\mathbf{W}(n;\underline{1}))[[t]]\otimes
I_t.
\end{split}
\end{equation*}

Thereby, we prove that the ideal $I_t$ is also a coideal of the Hopf
algebra $ U(\mathbf{W}(n;\underline{1}))[[t]]$.

\smallskip
(II) \ By Lemmas 3.9, 3.2, 3.4 \& 3.10, we have
\begin{equation*}
\begin{split}
S((x^{(\alpha)}D_i)^p)
&=-(1{-}e(k)t)^{-p(\alpha_k-\delta_{ik})}(1{-}e(m)t)^{-p(\alpha_m-\delta_{im})}\\
& \quad \cdot \Bigl( \sum\limits_{n,\ell=0}^{\infty}
d_m^{(n)}d_k^{(\ell)}((x^{(\alpha)}D_i)^p)\cdot h(k)_1^{\lg
\ell\rg}h(m)_1^{\lg
n\rg}t^{n+\ell}\Bigr)\\
&\equiv -\sum\limits_{n,\ell=0}^{p-1}
d_m^{(n)}d_k^{(\ell)}((x^{(\alpha)}D_i)^p)\cdot h(k)_1^{\lg
\ell\rg}h(m)_1^{\lg n\rg}t^{n+\ell} \quad
(\text{mod }(p, I))\\
&=-(x^{(\alpha)}D_i)^p+\delta_{ik}\delta_{\alpha,\epsilon_i}e(i)\cdot
h(k)_1^{\lg 1\rg} t+\delta_{im}\delta_{\alpha,\epsilon_i}e(i)\cdot
h(m)_1^{\lg 1\rg} t.
\end{split}\tag{54}
\end{equation*}
Hence, when $\alpha\ne\epsilon_i$, we get
$$S\bigl((x^{(\al)}D_i)^p\bigr)\equiv-(x^{(\al)}D_i)^p\equiv 0\quad (\text{\rm mod} \, I_t).$$
When $\al=\epsilon_i$, by Lemmas 3.4 and 3.9, (48) reads as
$S(x^{(\epsilon_i)}D_i)=-x^{(\epsilon_i)}D_i+\delta_{ik}e(i)\cdot
h(k)_1^{\langle 1\rangle}t+\delta_{im}e(i)\cdot h(m)_1^{\langle
1\rangle}t$. Combining with (54), we obtain
$$S\bigl((x^{(\epsilon_i)}D_i)^p-x^{(\epsilon_i)}D_i\bigr)\equiv-\bigl((x^{(\epsilon_i)}D_i)^p-x^{(\epsilon_i)}D_i\bigr)\equiv
0\quad (\text{\rm mod} \, I_t).$$

Thereby, we verify that $I_t$ is preserved by the antipode $S$ of
the quantization $U(\mathbf{W}(n;\underline{1}))[[t]]$ given in
Lemma 3.8.

The proof is complete.
\end{proof}

\begin{remark}
Theorem 3.7 gives $2^n{-}1$ new Hopf algebra structures of
$\mathcal{K}$-dimension $p^{np^n}$ over the same restricted
universal enveloping algebra
$\mathbf{u}(\mathbf{W}(n;\underline{1}))$ under the assumption that
$\mathcal{K}$ is algebraically closed and that $t$ is specialized at
a root of the $p$-polynomial $t^p-q\,t$ with $q\in\mathcal{K}^*$.
\end{remark}

\vskip30pt \centerline{\bf ACKNOWLEDGMENT}

\vskip15pt NH would like to express his thanks to the ICTP for its
support when he visited the ICTP Mathematics Group from March to
August, 2006. Authors are indebted to David Hernandez, Leonard L.
Scott for their kind help in improving our English writing, to the
referee for his/her useful comments.

\bigskip
\bigskip

\bibliographystyle{amsalpha}

\begin{thebibliography}{A}

\medskip
\bibitem {CP}V. Chari and A. Pressley, \textit{A Guide to Quantum Groups}, Cambridge
University Press, Cambridge, 1995.

\bibitem{DZ}D. Dokovic and K. Zhao, \textit{Derivations, isomorphisms and
second cohomology of generalized-Witt algebras}, Trans. Amer.
Math. Soc. \textbf{350} (1998), 643--664.

\bibitem {D}V.G. Drinfel'd, \textit{Quantum groups}, Proceedings ICM
(berkeley 1986) \textbf{1} (1987), AMS, 798--820.

\bibitem {ES} P. Etingof and O. Schiffmann, \textit{Lectures on Quantum Groups}, 2nd,
International Press, USA, 2002.

\bibitem{AJ}A. Giaquinto and J. Zhang, \textit{Bialgebra action, twists and
universal deformation formulas}, J. Pure Appl. Algebra,
\textbf{128} (2) (1998), 133--151.

\bibitem{CG}C. Grunspan, \textit{Quantizations of the Witt algebra and of simple
Lie algebras in characteristic $p$}, J. Algebra, \textbf{280}
(2004), 145--161.

\bibitem{M} W. Michaelis, \textit{A class of infinite-dimensional Lie bialgebras containing
the Virasoro algebras}, Adv. Math., \textbf{107} (1994), 365--392.

\bibitem{NT} S.H. Ng and E.J. Taft, \textit{Classification of the Lie bialgebra structures
on the Witt and Virasoro algebras}, J. Pure and Appl. Algebra,
\textbf{151} (2000), 67--88.

\bibitem{DR}D.E. Radford, \textit{Operators on Hopf algebras}, Amer. J. Math.
\textbf{99} (1977), 139--158.

\bibitem{GY}G. Song and Y. Su, \textit{Lie
bialgeras of generalized-Witt type}, arXiv.Math: QA/0504168, Science
in China, Ser. A---Math. \textbf{49} (4) (2006), 533--544.


\bibitem{H}H. Strade, \textit{Simple Lie Algebras over Fields of Positive Characteristic,
I. Structure Theory}, de Gruyter Expositions in Mathematics, {\bf
38}, Walter de Gruyter, 2004.

\bibitem{HR}H. Strade and R. Farnsteiner, \textit{Modular Lie Algebras and
Their Representations}, Monogr. Textbooks, Pure Appl. Math. {\bf
116}, Marcel Dekker, 1988.

\bibitem{T}E. Taft, \textit{Witt and Virasoro algebras as bialgebras}, J. Pure Appl.
Algebra \textbf{87} (3) (1993), 301--312.
\end{thebibliography}

\end{document}